\documentstyle[12pt]{article}

\begin{document}

\date{}

\title{\bf Quantum Computing and the Jones Polynomial}

\author{Louis H. Kauffman \\
  Department of Mathematics, Statistics and Computer Science \\
  University of Illinois at Chicago \\
  851 South Morgan Street\\
  Chicago, IL, 60607-7045}

 \maketitle
  
 \thispagestyle{empty}
 
 \subsection*{\centering Abstract}

{\em
 This paper is an exploration of relationships between the Jones polynomial and quantum computing.
 We discuss the structure of the Jones polynomial in relation to representations of the Temperley Lieb algebra, and 
 give an example of a unitary representation of the braid group. We discuss the evaluation of the
 polynomial as a generalized quantum amplitude and show how the braiding part of the evaluation can be construed as a
 quantum computation when the braiding representation is unitary. The question of an efficient quantum 
 algorithm for computing the whole polynomial remains open.
}
 
\section{Introduction}
 This paper is an exploration of issues interrelating the Jones polynomial \cite{JO86} and quantum computing.
In  section 2 of the paper we review the formalism of Dirac brackets and some of the quantum physics associated
with this formalism. The section ends with a brief description of the concept of quantum computer that we shall
use in this paper. In section 3 we discuss the Jones and Temperley Lieb algebras and how they can be used to
produce representations of the Artin Braid group. While most of these representations are not unitary, we show
how to construct non-trivial unitary representations of the three-strand braid group by considering the
structure of two projectors. It turns out that two elementary projectors naturally generate a Temperley Lieb
algebra. This provides a way to make certain unitary representations and to motivate the construction of both
the Alexander and the Jones polynomial. In regard to the Alexander polynomial, we end this section with a
representation of the Artin Braid Group, constructed using projectors, that is equivalent to the classical
Burau representation. In section 4 we construct the bracket polynomial model for the Jones polynomial and
relate its structure to the representations discussed in the previous section. Section 5 shows how to
reformulate the bracket state sum in terms of discrete quantum amplitudes. This sets the stage for our
proposal, explained in section 6, for regarding knot invariants as quantum computers. This proposal needs
unitary braiding (a special condition) and the results of the computer are probabilistic. Nevertheless, I
believe that this model deserves consideration. The dialogue between topology and quantum computing is just
beginning. 
\vspace{3mm}

\noindent {\bf Acknowledgement.} Research on this paper was supported by National Science Foundation Grant DMS 9802859.

\section{Dirac Brackets}

 We begin with a discussion of  Dirac's notation,  $<b|a>$,  \cite{D58}.  In this notation
 $<a|$   and  $|b>$  are covectors and vectors respectively.    $<b|a>$  is the
evaluation of  $|a>$  by $<b|$,  hence it is a scalar, and in ordinary quantum
mechanics it is a complex number.   One can think of this as the amplitude for
the state to begin in Ò$a$Ó  and end in Ò$b$Ó.   That is, there is a process that can
mediate a transition from state  $a$ to state $b.$  Except for the fact that
amplitudes are complex valued, they obey the usual laws of probability.  This
means that if the process can be factored into a set of all possible intermediate
states  $c_{1}$, $c_{2}$, ..., $c_{n}$ ,  then the amplitude  for  $a \longrightarrow b$  is the
sum of the amplitudes for $a \longrightarrow c_{i} \longrightarrow b$.   Meanwhile, the
amplitude for $a \longrightarrow c_{i} \longrightarrow b$  is the product of the amplitudes of
the two subconfigurations $a \longrightarrow c_{i}$ and $c_{i} \longrightarrow b.$  Formally we
have

$$<b|a>  = \Sigma_{i} <b|c_{i}><c_{i}|a>$$

\noindent where the summation is over all the intermediate states  $i=1$, ..., $n.$

\noindent In general, the amplitude for mutually disjoint processes is the sum  of the
amplitudes of the individual processes.  The amplitude for a configuration of
disjoint processes is the product  of their individual amplitudes.
\vspace{3mm}

Dirac's division of the amplitudes into bras  $<b|$   and  kets   $|a>$  is done
mathematically by taking a vector space  $V$ (a Hilbert space, but it can be finite
dimensional) for the kets:  $|a>$  belongs to $V.$   The dual space $V^{*}$  is the home
of the bras.   Thus  $<b|$ belongs to $V^{*}$  so that   $<b|$ is a linear mapping 
$<b|:V \longrightarrow C$  where  $C$   denotes the complex numbers. We restore symmetry to the
definition by realising that an element of a vector space $V$ can be regarded as a
mapping from the complex numbers to $V.$  Given  $|a>: C \longrightarrow V$,  the
corresponding element of $V$ is the image of $1$ (in $C$)  under this mapping.  In
other words,  $|a>(1)$  is a member of $V.$   Now we have   $|a> :C \longrightarrow V$   and $<b|
: V \longrightarrow C.$ The composition  $<b|\circ|a> = <b|a> : C \longrightarrow C$ is regarded
as an element of $C$ by taking  the specific value $<b|a>(1).$   The complex numbers are
regarded as the ÒvacuumÓ, and the entire amplitude  $<b|a>$  is a Òvacuum to
vacuumÓ  amplitude for a process that includes the creation of the state $a$, its
transition to $b$, and the annihilation of $b$ to the vacuum once more.
\vspace{3mm}

Dirac notation has a life of its own. Let $$P =  |y><x|.$$ \noindent Let  $$<x| |y> = <x|y>.$$
\noindent Then  $$PP =  |y><x| |y><x|   = |y> <x|y> <x|   =  <x|y> P.$$ Up to a scalar
multiple,  $P$  is a projection operator. That is, if we let $$Q= P/<x|y>,$$ \noindent then $$QQ =
PP/<x|y><x|y> = <x|y>P/<x|y><x|y> = P/<x|y> = Q.$$ Thus $QQ=Q.$  In this language,
the completeness of intermediate states becomes the statement that a certain sum
of projections is equal to the identity: Suppose that $\Sigma_{i}  |c_{i}><c_{i}|   =  1$  
(summing over $i$) with  $<c_{i}|c_{i}>=1$ for each $i.$   Then

$$<b|a> =  <b| |a> =  <b|\Sigma_{i} |c_{i}><c_{i}|  |a> =  \Sigma_{i} <b| |c_{i}><c_{i}| |a>$$

$$<b|a> =  \Sigma_{i} <b|c_{i}><c_{i}|a>$$

\noindent Iterating this principle of expansion over a complete set of states leads to the
most primitive form of the Feynman integral \cite{FEY65}.   Imagine that the initial
and final states  $a$  and $b$  are points on the vertical lines  $x=0$  and $x=n+1$
respectively in the $x-y$ plane, and that  $(c(k)_{i(k)} , k)$  is a given point on the
line  $x=k$  for $0<i(k)<m.$  Suppose that   the sum   of projectors for each
intermediate state is complete.  That is, we assume that following sum is equal
to one, for each $k$  from $1$ to $n-1:$

$$|c(k)_{1}><c(k)_{1}| + ... + |c(k)_{m}><c(k)_{m}|  = 1.$$

\noindent Applying the completeness iteratively, we obtain the following expression for the
amplitude   $<b|a>:$

$$<b|a>  = \Sigma \Sigma \Sigma ... \Sigma <b|c(1)_{i(1)}><c(1)_{i(1)}|c(2)_{i(2)}> ... <c(n)_{i(n)}|a>$$

\noindent where the sum is taken over all  $i(k)$ ranging between  $1$  and $m$,  and  $k$ 
ranging
 between  $1$  and  $n.$   Each term in this sum can be construed as a combinatorial
path from  $a$  to  $b$  in the two dimensional space of the $x-y$ plane.     Thus the
amplitude for going from  $a$ to $b$   is seen as a summation of contributions from
all the ÒpathsÓ  connecting  $a$ to $b.$ See Figure 1.
\vspace{3mm}

{\tt    \setlength{\unitlength}{0.92pt}
\begin{picture}(432,308)
\thinlines    \put(397,279){\makebox(30,28){$C_{n}$}}
              \put(402,1){\makebox(29,27){$C_{1}$}}
              \put(303,150){\makebox(37,43){$B$}}
              \put(1,154){\makebox(41,37){$A$}}
              \put(359,274){\line(1,0){32}}
              \put(359,233){\line(1,0){32}}
              \put(359,193){\line(1,0){32}}
              \put(358,152){\line(1,0){32}}
              \put(359,113){\line(1,0){32}}
              \put(359,73){\line(1,0){32}}
              \put(359,32){\line(1,0){32}}
              \put(358,13){\line(1,0){32}}
              \put(358,52){\line(1,0){32}}
              \put(358,93){\line(1,0){32}}
              \put(358,133){\line(1,0){32}}
              \put(357,173){\line(1,0){32}}
              \put(357,213){\line(1,0){32}}
              \put(359,253){\line(1,0){32}}
              \put(261,347){\line(0,0){0}}
              \put(358,294){\line(1,0){32}}
              \put(373,293){\line(0,-1){280}}
              \put(173,29){\line(5,6){118}}
              \put(52,173){\line(5,-6){120}}
              \put(51,173){\line(0,-1){1}}
              \put(173,75){\line(6,5){116}}
              \put(51,174){\line(5,-4){122}}
              \put(173,112){\line(2,1){116}}
              \put(53,173){\line(2,-1){120}}
              \put(173,156){\line(6,1){119}}
              \put(174,196){\line(-5,-1){122}}
              \put(274,173){\line(0,0){0}}
              \put(53,175){\line(6,-1){121}}
              \put(174,196){\line(5,-1){115}}
              \put(173,235){\line(2,-1){119}}
              \put(51,173){\line(2,1){122}}
              \put(173,274){\line(6,-5){119}}
              \put(52,173){\line(6,5){121}}
              \put(172,10){\line(3,4){119}}
              \put(51,174){\line(3,-4){121}}
              \put(173,54){\line(1,1){118}}
              \put(52,174){\line(1,-1){120}}
              \put(173,93){\line(3,2){119}}
              \put(52,174){\line(3,-2){120}}
              \put(173,134){\line(3,1){119}}
              \put(52,174){\line(3,-1){120}}
              \put(171,296){\line(1,-1){121}}
              \put(52,174){\line(1,1){120}}
              \put(173,254){\line(3,-2){120}}
              \put(52,174){\line(3,2){120}}
              \put(172,214){\line(3,-1){121}}
              \put(52,174){\line(3,1){120}}
              \put(52,174){\line(1,0){242}}
              \put(293,172){\circle*{20}}
              \put(52,173){\circle*{20}}
\end{picture}}

\begin{center}
{\bf Figure 1 - Intermediates}
\end{center}
 \vspace{3mm}

Feynman used this description to produce his famous path integral expression for 
amplitudes in  quantum mechanics.  His path integral takes the form

$$\int dP exp(iS)$$

\noindent where  $i$  is the square root of minus one,  the integral is taken over all paths
from  point  $a$ to point $b$,  and   $S$  is the action  for a particle to travel from
 $a$ to $b$ along a given path.   For the quantum mechanics associated with a
classical (Newtonian) particle  the action  $S$  is given by  the integral along
the given path from a to b  of the difference $T-V$  where $T$ is the classical
kinetic energy and $V$ is the classical potential energy of the particle.
\vspace{3mm}

\subsection{What is a Quantum Computer?}
We are now in a position to explain the definition of quantum computer that will be used in this paper.
Let $H$ be a given finite dimensional vector space over the complex numbers $C.$ Let $\{ W_{0}, W_{1},..., W_{n} \}$ be an
orthonormal basis for $H$ so that with $|i> := |W_{i}>$ denoting $W_{i}$ and $<i|$ denoting the conjugate transpose of $|i>$,
we have
$$<i|j> = \delta_{ij}$$
\noindent where $\delta_{ij}$ denotes the Kronecker delta (equal to one when its indices are equal to one another, and equal
to zero otherwise). Given a vector $v$ in $H$ let $|v|^{2} := <v|v>.$ Note that $<i|v$ is the $i$-th coordinate of $v.$ 
\vspace{3mm}

\noindent An {\em observation of $v$} returns one of the coordinates $|i>$
of $v$ with probability $|<i|v|^{2}.$ This model of observation is a simple instance of the situation with a quantum
mechanical system that is in a mixed state until it is observed. The result of observation is to put the system into one of
the basis states. 
\vspace{3mm}

When the dimension of the space $H$ is two ($n=1$), a vector in the space is called a {\em qubit}. A qubit represents one
quantum of binary information. On observation, one obtains either the ket $|0>$ or the ket $|1>$. This constitutes the 
binary distinction that is inherent in a qubit.  Note however that the information obtained is probabilistic.  If the qubit is
$$\psi = \alpha |0> + \beta \ |1>,$$ \noindent then the ket $|0>$ is observed with probability $|\alpha|^{2}$, and the ket
$|1>$ is observed with probability $|\beta|^{2}.$  In speaking of an idealized quantum computer, we do not specify the nature
of measurement process beyond these probability postulates.
\vspace{3mm}
 
In the case of general dimension $n$ of the space $H$, we will call the vectors in $H$
{\em qunits}. It is quite common to use spaces $H$ that are tensor products of two-dimensional spaces (so that all computations 
are expressed in terms of qubits) but this is not neccessary in principle. One can start with a given space, and later work out
factorizations into qubit transformations.
\vspace{3mm}

A {\em quantum computation} consists in the application of a unitary
transformation $U$ to an initial qunit $\psi = a_{1}|1> + ... + a_{n}|n>$  with $|\psi|^{2}=1$, plus an
observation  of
$U\psi.$ An observation of $U\psi$ returns the ket $|i>$ with probability $|U\psi|^{2}$. In particular, if we start the computer
in the state $|i>$, then the probability that it will return the state $|j>$ is $|<j|U|i>|^{2}.$

\vspace{3mm} It is the neccessity for writing a given computation in terms of unitary trasformations, and the probabilistic
nature of the result that characterizes quantum computation. Such computation could be carried out by an idealized quantum
mechanical system. It is hoped that such systems can be physically realized. 
\vspace{3mm}

\section{Braiding, Projectors and the Temperley Lieb Algebra}

The Jones polynomial is one 
of the great mathematical breakthroughs of the twentieth century, and like many such breakthroughs 
it appears basically simple in retrospect. I will tell two stories in this section. The first story
is a capsule summary of how Jones discovered the polynomial by way of an apparently strange 
algebraic structure that first appeared in his research on von Neumann algebras, and then was pointed
out to be an algebra known to experts in the Potts model in statistical mechanics.   The second story shows that the essential
algebra for the needed representation of the braid group is present in the algebra generated by any two simple projectors (see
below for the definitions of these terms) and that it is graphically illustrated by the Dirac bra-ket notation for these
operators.
\vspace{3mm}

Jones was studying the inclusion of one von Neuman algebra $N$ in another one $M$. In this 
context there is a projection $e_{1}:M \longrightarrow N$ so that the restriction of $e$ to $N$ is the 
identity mapping, and so that $e_{1}^{2} = e_{1}$. In his context the algebra M could be extended to 
include this projector to an algebra $M_{1} = M \bigcup \{e_{1} \}.$ Then we have
$$N \subset M \subset M_{1}$$
\noindent and the construction can be continued inductively to produce
$$N \subset M \subset M_{1} \subset M_{2} \subset M_{3} \subset ...$$
\noindent and an algebra of projectors
$$e_{1}, e_{2}, e_{3},...$$
\noindent such that 
$$e_{i}^{2} = e_{i}, i = 1,2,3,...$$
$$e_{i}e_{i \pm 1}e_{i} = \kappa e_{i}, i = 2,3,...$$
$$e_{i}e_{j} = e_{j}e_{i}, \hspace{1mm} |i-j|>1.$$
\vspace{3mm}

\noindent We will call an algebra that can be expressed with generators and relations as above a {\em Jones algebra}. 
$J_{\infty}$ will denote a
Jones algebra on infinitely many generators as above. $J_{n}$ will denote the Jones algebra generated by an identity
element $1$ and generators $e_{1}, ..., e_{n-1}.$

It was pointed out that the relations
$$e_{i}e_{i \pm 1}e_{i} = \kappa e_{i}, i = 2,3,...$$
$$e_{i}e_{j} = e_{j}e_{i} , |i-j|>1$$
\noindent look suspiciously like the basic braiding relations in the Artin Braid group which read
$$\sigma_{i}\sigma_{i \pm 1}\sigma_{i} =\sigma_{i \pm 1}\sigma_{i}\sigma_{i \pm 1} , i = 2,3,...$$
$$\sigma_{i}\sigma_{j} = \sigma_{j}\sigma_{i} , |i-j|>1.$$
\vspace{3mm}

{\tt    \setlength{\unitlength}{0.92pt}
\begin{picture}(405,307)
\thicklines   \put(163,1){\makebox(81,56){$\sigma_{2}^{-1}$}}
              \put(1,3){\makebox(70,54){$\sigma_{1}^{-1}$}}
              \put(349,179){\makebox(40,39){$1$}}
              \put(183,175){\makebox(48,43){$\sigma_{2}$}}
              \put(12,176){\makebox(40,37){$\sigma_{1}$}}
              \put(163,144){\line(0,-1){80}}
              \put(243,103){\line(0,-1){40}}
              \put(203,103){\line(0,-1){40}}
              \put(82,144){\line(0,-1){80}}
              \put(43,105){\line(0,-1){41}}
              \put(4,104){\line(0,-1){41}}
              \put(243,143){\line(-1,-1){40}}
              \put(204,143){\line(1,-1){15}}
              \put(228,118){\line(1,-1){15}}
              \put(28,119){\line(1,-1){15}}
              \put(4,144){\line(1,-1){15}}
              \put(43,144){\line(-1,-1){40}}
              \put(402,304){\line(0,-1){80}}
              \put(363,304){\line(0,-1){81}}
              \put(322,304){\line(0,-1){81}}
              \put(243,265){\line(0,-1){42}}
              \put(204,263){\line(0,-1){39}}
              \put(163,304){\line(0,-1){79}}
              \put(82,304){\line(0,-1){80}}
              \put(43,264){\line(0,-1){40}}
              \put(4,263){\line(0,-1){38}}
              \put(205,304){\line(1,-1){39}}
              \put(204,264){\line(1,1){15}}
              \put(243,302){\line(-1,-1){13}}
              \put(42,303){\line(-1,-1){13}}
              \put(3,265){\line(1,1){15}}
              \put(4,305){\line(1,-1){39}}
\end{picture}}

\begin{center}
{\bf Figure 2 - Braid Group Generators} 
\end{center}
\vspace{3mm} 

This pattern led Jones to first construct a representation of the Artin Braid Group to his algebra, and then to discover
an invariant of knots and links that is related to this representation. Figure 2 illustrates the generators of the braid
group. The second and third Reidemeister moves shown in Figure 6 illustrate the braiding relations except for commutativity
of distant generators.
\vspace{3mm}

The representation that Jones discovered is a linear one in the form of 
$$\rho: B_{\infty} \longrightarrow J_{\infty}$$
\noindent where 
$$\rho(\sigma_{i}) = \alpha 1 + \beta e_{i}$$
\noindent for appropriate constants $\alpha$ and $\beta.$ We will elaborate on this representation shortly.
Here $J_{\infty}$ denotes the algebra generated by the $e_{i}$ for $i=1,2,3,...$.
It seems an amazing coincidence that a representation algebra for the Artin Braid group would appear in a context
that seems so far away from this structure. The complex source of Jones' algebra makes this connection seem quite
mysterious, and the fact that this same algebra appears in statistical mechanics also seems mysterious. What is the
source of this apparent connection of the Artin Braid group with algebras and structures coming from quauntum physics?
\vspace{3mm}

\noindent {\bf Remark.} In the discussion to follow, we will use the bra and ket notations of Dirac and we will write
$<v|$ for $v^{t}$, using the notation $v^{t}$ for the transpose of a vector $v$.  It is to be understood that in the case of
a complex vector space, this is the conjugate transpose, but that in the generalizations that we use (over more general rings)
we will simply take the formal transpose without conjugation. Later we will construct real-valued representations of the 
Temperley-Lieb algebra, and there transpose will be the same as conjugate transpose.
\vspace{3mm}
 
For the purpose of this discussion it will be useful to define a {\em projector} to be a linear map $P:V
\longrightarrow V$ where $V$ is a vector space or a module over a ring $k$, and $P^{2}$ is a non-zero multiple of $P$. Shall
call a projector {\em simple} if, in a basis, it takes the form $P = vv^{t}$ where $v$ is a column vector and $v^{t}$ is its
transpose. Then $v^{t}v$ is the dot product of $v$ with itself and hence a scalar. Therefore
$$P^{2} = PP = vv^{t}vv^{t} = v[v^{t}v]v^{t}$$ 
$$= [v^{t}v]vv^{t} = [v^{t}v]P.$$
\noindent Because of the ubiquity of projectors in quantum physics, the physicist P.A.M. Dirac devised a beautiful notation for
this situation. Dirac would write $|v>$ for $v$ and $<v|$ for $v^{t}.$  He would write $$<v||w>=<v|w> = v^{t}w$$
\noindent for the dot product of two vectors in a given basis.
\vspace{3mm}

Then one can write $P$ in Dirac notation by the formula
$$P = |v><v|$$
\noindent and we have
$$P^{2} = PP = |v><v||v><v| = |v><v|v><v| $$
$$= <v|v>|v><v| = <v|v>P.$$

Now {\em consider the algebra generated by two simple projectors} $P=|v><v|$ and $Q=|w><w|.$ We have
$$P^{2} = <v|v>P,$$
$$Q^{2} = <w|w>Q$$
\noindent and
$$PQP = |v><v||w><w||v><v|$$
$$=|v><v|w><w|v><v|$$ 
$$=<v|w><w|v>|v><v|$$ 
$$= <v|w><w|v>P$$

\noindent while
$$QPQ = |w><w||v><v||w><w|$$
$$=|w><w|v><v|w><w|$$
$$=<w|v><v|w>|w><w|$$
$$= <w|v><v|w>Q$$ 
$$=<v|w><w|v>Q. $$
\noindent Thus, with $\lambda = <v|w><w|v>$ we have that 
$$PQP = \lambda P$$
$$QPQ = \lambda \, Q.$$
\noindent We can define $e = P/<v|v>$ and $f=Q/<w|w>$ and find
$$e^{2} = e$$
$$f^{2} = f$$
$$efe = \kappa e$$
$$fef = \kappa f$$
\noindent where $\kappa = \lambda/(<v|v><w|w>).$ In this way we see that {\em any two simple projectors generate a Jones
algebra of type $J_{2}$.} In this sense the appearance of such algebras is quite natural.  The relationship with braiding
remains as remarkable as ever.
\vspace{3mm}

In order to see how these representations work, it is useful to discuss the combinatorics of these algebras a bit
further. The {\em Temperley Lieb algebra} $TL_{n}$ \cite{KA87} is an algebra over a commutative ring $k$ with
generators $\{ 1, U_{1},U_{2}, ... ,U_{n-1} \}$ and relations 

$$U_{i}^{2} = \delta U_{i},$$

$$U_{i}U_{i \pm 1}U_{i} = U_{i},$$

$$U_{i}U_{j} = U_{j}U_{i}, |i-j|>1,$$

\noindent where $\delta$ is a chosen element of the ring $k$. These equations give the
multiplicative structure of the algebra. The algebra is a free module over the ring $k$ with basis
the equivalence classes of these products modulo the given relations.
\vspace{3mm}

We will make the ground ring specific in the examples to follow. It is clear that the concepts
of Temperley Lieb algebra and Jones algebra are interchangeable. Given a Jones algebra $J_{\infty}$, with
$e_{i}e_{i \pm 1}e_{i} = \kappa e_{i}$, let $\delta = 1/\sqrt{\kappa}$ (assuming that this square root exists in the ground
ring $k$. Then let $U_{i} = \delta e_{i}$ and we find that $U_{i}^{2} = \delta U_{i}$ with 
$$U_{i}U_{i \pm 1}U_{i} = (1/ \sqrt{\kappa})^{3}e_{i}e_{i \pm 1}e_{i}$$
$$= (1/\sqrt{\kappa})^{3} \kappa e_{i} = (1/\sqrt{\kappa}) e_{i} = U_{i},$$
\noindent converting the Jones algebra to a Temperley Lieb algebra.
\vspace{3mm}

It is useful to see the bare bones of the algebra of two projectors. For this purpose, lets write
$$P= ><$$
\noindent and 
$$Q= ][.$$
\noindent Then 
$$PP = >< \, >< = <> \, >< = <> P$$
\noindent and 
$$QQ = ][ \, ][ = [] Q$$
\noindent while
$$PQP = >< \, ][ \, >< = <] \, [> P$$
$$QPQ = ][ \, >< \, ][ = [> \, <] Q = <] \,[> Q.$$
\noindent To see how the representation of the braid group is constructed, lets assume that the scalars
$<]$ and $[>$ are both equal to $1$ and that $\delta = <> = []$. Then $P$ and $Q$ form a two-generator Temperley Lieb
algebra $TL_{3}$. We will illustrate how to represent the three strand Artin braid group $B_{3}$ to $TL_{2}$. 
\vspace{3mm}

It is useful to use the iconic symbol $><$ for a projector and to choose another iconic symbol  \mbox{\large $\asymp$}
for the identity operator in the algebra. With these choices we have
$$\mbox{\large $\asymp$}\mbox{\large $\asymp$} \,\, = \,\, \mbox{\large $\asymp$}$$
$$\mbox{\large $\asymp$} >< \,\, = \,\, >< \mbox{\large $\asymp$} \,\,= \,\, ><$$
$$\mbox{\large $\asymp$} ][ \,\, = \,\, ][ \mbox{\large $\asymp$} \,\, = \,\, ][$$
\vspace{3mm}

\noindent We define the representation $\rho:B_{3} \longrightarrow TL_{3}$ on the generators $\sigma_{1} = \sigma$ and 
$\sigma_{2} = \tau$ of the three strand braid group, whose relations are $\sigma \tau \sigma = \tau \sigma \tau$ plus the
invertibility of the generators. We define
$$\rho(\sigma) = A 1 + B P = A \mbox{\large $\asymp$} + B ><$$
$$\rho(\sigma^{-1}) = B 1 + A P =  B \mbox{\large $\asymp$} + A ><$$
\noindent and
$$\rho(\tau) = A 1 + B Q= A \mbox{\large $\asymp$} + B ][$$
$$\rho(\tau^{-1}) = B 1 + A Q= B \mbox{\large $\asymp$} + A ][.$$
\noindent where $A$ and $B$ are commuting indeterminates.
\vspace{3mm}

With these definitions, we have 
$$\rho(\sigma) = A \mbox{\large $\asymp$} + B >< $$
$$\rho(\sigma^{-1}) = B \mbox{\large $\asymp$} + A >< .$$
\noindent Thus
$$ \mbox{\large $\asymp$} =( A \mbox{\large $\asymp$} + B ><)(B \mbox{\large $\asymp$} + A ><)$$ 
$$= AB \mbox{\large $\asymp$}\mbox{\large $\asymp$} + A^{2}\mbox{\large $\asymp$}>< + B^{2}><  \mbox{\large $\asymp$}
+ AB ><><$$
$$= AB \mbox{\large $\asymp$}  + A^{2}>< + B^{2}>< + AB \delta ><$$
$$\mbox{\large $\asymp$} = AB \mbox{\large $\asymp$}  + (A^{2} + B^{2} + AB \delta) ><$$
\noindent Consequently, we will have $1 = \rho(\sigma)\rho(\sigma^{-1})$ if we take 
$B=A^{-1}$
\noindent and
$\delta = -A^{2} - A^{-2}.$
\noindent We shall take these values from now on so that
$$\rho(\sigma) = A \mbox{\large $\asymp$} + A^{-1} >< = A 1 + A^{-1} P$$
\noindent and 
$$\rho(\tau) = A \mbox{\large $\asymp$} + A^{-1} ][ = A 1 + A^{-1} Q. $$
\vspace{3mm}

With these specializations of $A$ and $B$, it is easy to verify that $\rho$ is a representation of the Artin Braid Group.
Note that $P^{2}= \delta P$, $Q^{2}=\delta Q$ and $PQP = P.$ 
$$\rho(\sigma) \rho(\tau) \rho(\sigma) = (A + A^{-1}P)(A + A^{-1}Q)(A + A^{-1}P)$$
$$=(A^{2} + Q + P + A^{-2}PQ)(A + A^{-1}P)$$
$$=A^{3} + AQ + AP + A^{-1}PQ + AP + A^{-1}QP + A^{-1}P^{2} + A^{-3}PQP$$
$$=A^{3} + AQ + AP + A^{-1}PQ + AP + A^{-1}QP + A^{-1}\delta P + A^{-3}P$$
$$=A^{3}  + (2A + + A^{-1}(-A^{2}-A^{-2}) + A^{-3})P + AQ + A^{-1}(PQ + QP)$$
$$=A^{3}  + AP + AQ + A^{-1}(PQ + QP)$$
$$ \rho(\sigma) \rho(\tau) \rho(\sigma) = A^{3}  + A(P + Q) + A^{-1}(PQ + QP)$$
\noindent Since this last expression is symmetric in $P$ and $Q$, we conclude that 
$$\rho(\sigma) \rho(\tau) \rho(\sigma) = \rho(\tau) \rho(\sigma) \rho(\tau).$$
\noindent Hence $\rho$ is a representation of the Artin
Braid Group.
\vspace{3mm}

This argument generalizes to yield a corresponding representation of the Artin Braid Group $B_{n}$ to the Temperley Lieb
algebra $TL_{n}$ for each $n=2,3,....$ We will discuss the structure of these representations below. In the next
section we show how the Jones polynomial can be constructed by a state summation model. This model can be also be viewed as
a generalisation of the above representation of the Temperley Lieb algebra. 
\vspace{3mm}

The very close relationship between elementary quantum mechanics and topology is very well
illustrated by the structure and representations of the Temperley Lieb algebra. 
\vspace{3mm}

{\tt    \setlength{\unitlength}{0.92pt}
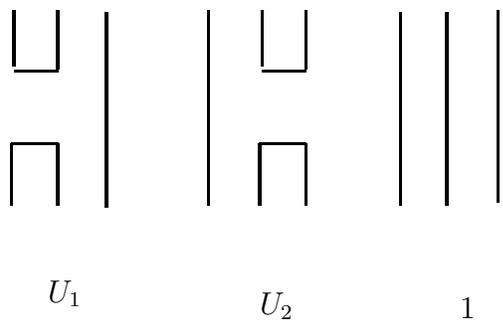
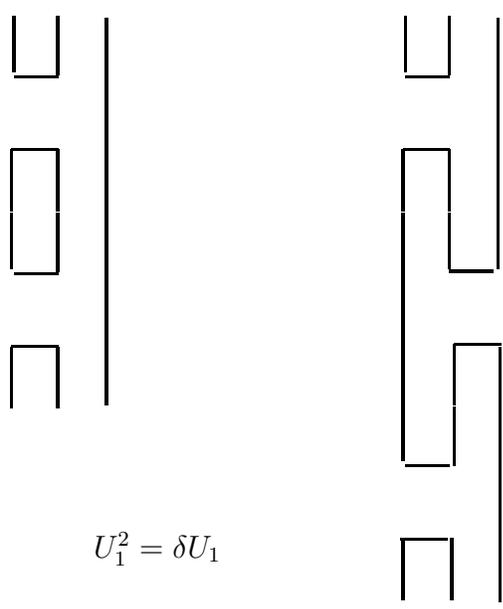
\begin{picture}(273,523)
\thicklines   \put(110,1){\makebox(162,80){$U_{1}U_{2}U_{1}=U_{1}$}}
              \put(1,93){\makebox(180,100){$U_{1}^{2}=\delta U_{1}$}}
              \put(191,374){\makebox(55,51){$1$}}
              \put(111,375){\makebox(58,51){$U_{2}$}}
              \put(22,382){\makebox(62,47){$U_{1}$}}
              \put(192,202){\line(0,-1){23}}
              \put(192,147){\line(0,-1){25}}
              \put(232,201){\line(0,-1){80}}
              \put(213,201){\line(0,-1){24}}
              \put(212,122){\line(0,1){25}}
              \put(193,177){\line(1,0){18}}
              \put(191,147){\line(1,0){19}}
              \put(213,227){\line(1,0){19}}
              \put(211,257){\line(1,0){18}}
              \put(232,201){\line(0,1){25}}
              \put(231,282){\line(0,-1){24}}
              \put(192,281){\line(0,-1){80}}
              \put(213,227){\line(0,-1){25}}
              \put(211,281){\line(0,-1){23}}
              \put(193,362){\line(0,-1){23}}
              \put(192,307){\line(0,-1){25}}
              \put(231,361){\line(0,-1){80}}
              \put(211,362){\line(0,-1){24}}
              \put(211,282){\line(0,1){25}}
              \put(193,337){\line(1,0){18}}
              \put(192,307){\line(1,0){19}}
              \put(231,522){\line(0,-1){79}}
              \put(210,521){\line(0,-1){79}}
              \put(191,521){\line(0,-1){79}}
              \put(134,522){\line(0,-1){23}}
              \put(133,467){\line(0,-1){25}}
              \put(112,522){\line(0,-1){80}}
              \put(152,522){\line(0,-1){24}}
              \put(152,442){\line(0,1){25}}
              \put(134,497){\line(1,0){18}}
              \put(133,467){\line(1,0){19}}
              \put(31,281){\line(0,-1){23}}
              \put(31,226){\line(0,-1){25}}
              \put(70,282){\line(0,-1){80}}
              \put(50,281){\line(0,-1){24}}
              \put(50,201){\line(0,1){25}}
              \put(32,256){\line(1,0){18}}
              \put(31,226){\line(1,0){19}}
              \put(32,362){\line(0,-1){23}}
              \put(31,307){\line(0,-1){25}}
              \put(70,361){\line(0,-1){80}}
              \put(50,362){\line(0,-1){24}}
              \put(50,282){\line(0,1){25}}
              \put(32,337){\line(1,0){18}}
              \put(31,307){\line(1,0){19}}
              \put(31,467){\line(1,0){19}}
              \put(32,497){\line(1,0){18}}
              \put(50,442){\line(0,1){25}}
              \put(50,522){\line(0,-1){24}}
              \put(70,521){\line(0,-1){80}}
              \put(31,467){\line(0,-1){25}}
              \put(32,522){\line(0,-1){23}}
\end{picture}}

\begin{center}
{\bf Figure 3 - Diagrammatic Temperley Lieb Algebra} 
\end{center}
\vspace{3mm}

Figure 3 illustrates a diagrammatic interpretation of the Temperley Lieb algebra. In this interpretation,
the multiplicative generators of the module are collections of strands connecting $n$ top points and $n$ bottom points. Top
points can be connected either to top or to bottom points. Bottom points can be connected to either bottom or to
top points. All connections are made in the plane with no overlapping lines and no lines going
above the top row of points or below the bottom row of points. Multiplication is accomplished by
connecting the bottom row of one configuration with the top row of another. In Figure 3 we have
illustrated the types of special configurations that correspond to the $U_{i}$, and we have shown
that $\delta$ is interpreted as a closed loop.
\vspace{3mm}

{\tt    \setlength{\unitlength}{0.92pt}
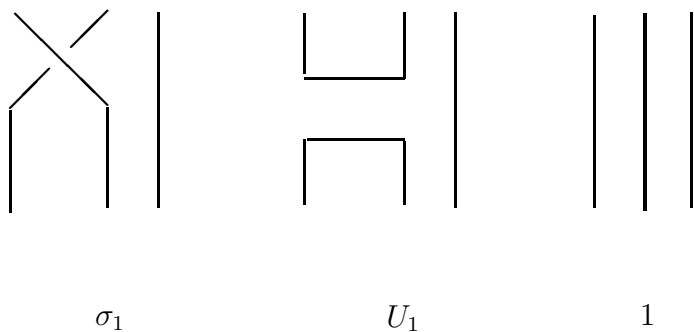
\begin{picture}(286,248)
\thicklines   \put(3,1){\makebox(282,84){$\rho(\sigma_{1})=AU_{1} + A^{-1}1$}}
              \put(245,91){\makebox(40,60){$1$}}
              \put(123,87){\makebox(83,65){$U_{1}$}}
              \put(4,87){\makebox(80,63){$\sigma_{1}$}}
              \put(283,245){\line(0,-1){80}}
              \put(264,245){\line(0,-1){81}}
              \put(243,244){\line(0,-1){79}}
              \put(186,245){\line(0,-1){80}}
              \put(165,192){\line(0,-1){26}}
              \put(125,193){\line(1,0){40}}
              \put(124,166){\line(0,1){27}}
              \put(165,218){\line(0,1){27}}
              \put(124,218){\line(1,0){41}}
              \put(124,245){\line(0,-1){25}}
              \put(64,245){\line(0,-1){80}}
              \put(43,206){\line(0,-1){41}}
              \put(3,205){\line(0,-1){42}}
              \put(28,231){\line(1,1){15}}
              \put(3,206){\line(1,1){16}}
              \put(5,245){\line(1,-1){38}}
\end{picture}}

\begin{center}
{\bf Figure 4 - Braid Group Representation} 
\end{center}
\vspace{3mm}

{\tt    \setlength{\unitlength}{0.92pt}
\begin{picture}(340,424)
\thinlines    \put(19,99){\makebox(100,57){$M^{ab}M_{ab}$}}
              \put(59,309){\makebox(59,74){$\delta^{b}_{a}$}}
              \put(238,265){\makebox(97,57){$M^{cd}$}}
              \put(252,358){\makebox(87,45){$M_{ab}$}}
              \put(6,1){\makebox(320,82){$M^{ai}M_{ib} = \delta^{a}_{b}$}}
\thicklines   \put(77,197){\makebox(20,20){$b$}}
              \put(1,199){\makebox(17,18){$a$}}
              \put(240,78){\makebox(15,21){$b$}}
              \put(197,157){\makebox(15,20){$i$}}
              \put(122,237){\makebox(14,20){$a$}}
              \put(237,316){\makebox(18,23){$d$}}
              \put(161,315){\makebox(19,22){$c$}}
              \put(237,357){\makebox(19,20){$b$}}
              \put(157,357){\makebox(20,19){$a$}}
              \put(66,209){\circle*{10}}
              \put(27,209){\circle*{10}}
              \put(228,91){\circle*{10}}
              \put(187,170){\circle*{10}}
              \put(148,249){\circle*{10}}
              \put(227,329){\circle*{10}}
              \put(188,329){\circle*{10}}
              \put(227,366){\circle*{10}}
              \put(187,366){\circle*{10}}
              \put(27,291){\circle*{10}}
              \put(28,411){\circle*{10}}
              \put(28,392){\makebox(31,31){$b$}}
              \put(27,274){\makebox(25,30){$a$}}
              \put(227,169){\line(0,-1){79}}
              \put(147,250){\line(0,-1){80}}
              \put(186,210){\line(1,0){41}}
              \put(227,209){\line(0,-1){39}}
              \put(186,210){\line(0,-1){41}}
              \put(147,170){\line(0,-1){40}}
              \put(148,130){\line(1,0){38}}
              \put(186,130){\line(0,1){40}}
              \put(26,249){\line(1,0){41}}
              \put(67,248){\line(0,-1){39}}
              \put(26,249){\line(0,-1){41}}
              \put(26,209){\line(0,-1){40}}
              \put(26,169){\line(1,0){38}}
              \put(65,170){\line(0,1){40}}
              \put(226,288){\line(0,1){40}}
              \put(188,288){\line(1,0){38}}
              \put(187,328){\line(0,-1){40}}
              \put(186,408){\line(0,-1){41}}
              \put(227,407){\line(0,-1){39}}
              \put(186,408){\line(1,0){41}}
              \put(26,410){\line(0,-1){120}}
\end{picture}}

\begin{center}
{\bf Figure 5 - Abstract Tensors} 
\end{center}
\vspace{3mm}

One way to make a matrix representation of the Temperley Lieb algebra (and a corresponding
representation of the braid group) is to use the matrix $M$ defined as follows

$$M =  \left[
\begin{array}{cc}
     0 & iA  \\
     -iA^{-1} & 0
\end{array}
\right]. $$

\noindent Note that $M^{2}=1$ where $1$ denotes the  ($ 2 \times 2$) identity matrix.
We will use $M$ with either upper or lower indices so that $M^{ab} = M_{ab}.$
$M$ will represent both the cup and the cap in the Temperley Lieb diagrams, with $M_{ab}$ representing the cap and 
$M^{ab}$ representing the cup.  If $U$ denotes a cup over a cap, then 

$$U^{ab}_{cd} = M^{ab}M_{cd}.$$

\noindent Note that 

$$(U^{2})^{ab}_{cd} = \Sigma_{ij} U^{ab}_{ij} U^{ij}_{cd} =$$

$$= \Sigma_{ij} M^{ab}M_{ij}M^{ij}M_{cd} = [\Sigma_{ij}M_{ij}M^{ij}]M^{ab}M_{cd}$$

$$= [\Sigma_{ij}M_{ij}M^{ij}]U^{ab}_{cd}.$$

\noindent Note that

$$\Sigma_{ij}M_{ij}M^{ij} = \Sigma_{ij}(M_{ij})^{2} = -A^{2}-A^{-2}.$$

\noindent Thus, letting $\delta = -A^{2}-A^{-2}$, we have

$$U^{2} = \delta U.$$

\noindent Then we take 
$U_{i}$	as a tensor product of identity matrices corresponding to the vertical lines in the diagram
for this element and one factor of $U$ for the placement of the cup-cap at the locations $i$ and
$i+1$. To see how this works to give the relation $U_{i}U_{i \pm 1}U_{i} = U_{i},$  we verify that 
$U_{1}U_{2}U_{1} = U_{1}$ in $TL_{3}.$ In the calculation to follow we will use the Einstein summation convention. Repeated
upper and lower indices are summed across the index set $\{1,2\}.$

$$U_{1} = U \otimes 1$$
\noindent and
$$U_{2} = 1 \otimes U$$
\noindent so that 

$$(U_{1})^{abc}_{def} = M^{ab}M_{de} \delta^{c}_{f}$$

\noindent and
$$ (U_{2})^{abc}_{def} = \delta^{a}_{d} M^{bc}M_{ef}.$$

\noindent Therefore

$$(U_{1}U_{2}U_{1})^{abc}_{def} = (U_{1})^{abc}_{ijk}(U_{2})^{ijk}_{rst}(U_{1})^{rst}_{def}$$

$$= ( M^{ab}M_{ij} \delta^{c}_{k})( \delta^{i}_{r}M^{jk}M_{st})(  M^{rs}M_{de} \delta^{t}_{f})$$

$$= M^{ab}(M_{rj}M^{jc})(M_{sf} M^{rs})M_{de} =  M^{ab}(\delta_{r}^{c})(\delta_{f}^{r})M_{de} = $$

$$= M^{ab}M_{de}\delta^{c}_{f} = (U_{1})^{abc}_{def}$$

\noindent Thus $$U_{1}U_{2}U_{1} = U_{1}.$$

\noindent This representation of the Temperley Lieb algebra is useful for knot theory and it is
conjectured to be a faithful representation. One may also conjecture that the corresponding braid
group representation is faithful.  \vspace{3mm}

\noindent {\bf Remark.} The reader should note that the diagrammatic interpretation of the Temperley Lieb algebra gives a
clear way to follow the index details of the calculation we have just performed. In the diagrams an index that is not on a
free end is summed over just as in the Einstein summation convention. An index at the end of a line is a free index and does
not receive summation. See Figure 5 for an illustration of the index algebra in relation to these diagrams. We will
generalize the diagrammatic algebra in section 5.
\vspace{3mm}

\subsection{Two Projectors and a Unitary Representation of the Three Strand Braid Group}
The Temperley Lieb representation of the braid group that we have described is not a unitary representation
except when $A^{2}=-1$, a value that is not of interest in the knot theory. In order to find
elementary unitary representations of the braid group, one has to go deeper.  
\vspace{3mm}

It is useful to think of the Temperley Lieb algebra as generated by projections 
$e_{i} = U_{i}/\delta$ so that $e_{i}^{2} = e_{i}$ and $e_{i}e_{i\pm 1}e_{i} = \tau e_{i}$ where
$\tau = \delta^{-2}$ and $e_{i}$ and $e_{j}$ commute for $|i-j|>1.$
\vspace{3mm} 

With this in mind, consider elementary projectors
$e = |A><A|$ and $f=|B><B|$. We assume that $<A|A> = <B|B> =1$ so that $e^{2} = e$ and $f^{2} =f.$
Now note that

$$efe = |A><A|B><B|A><A| = <A|B><B|A> e = \tau e$$

\noindent Thus $$efe=\tau e$$ 

\noindent where $\tau = <A|B><B|A>$.

This algebra of two projectors is the simplest instance of a representation of the Temperley Lieb
algebra.  In particular, this means that a representation of the three-strand braid group is
naturally associated with the algebra of two projectors, a simple toy model of quantum physics!
\vspace{3mm}

Quite specifically if we let $<A| = (a,b)$ and $|A> = (a,b)^{t}$ the transpose of this row vector,
then

$$e=|A><A| =  \left[
\begin{array}{cc}
     a^{2} & ab  \\
     ab & b^{2}
\end{array}
\right] $$ 

\noindent is a standard projector matrix when $a^{2} + b^{2} = 1.$  To obtain a specific
representation, let

$$e_{1} =  \left[
\begin{array}{cc}
     1 & 0  \\
     0 & 0
\end{array}
\right] $$ 

\noindent and 

$$e_{2} =  \left[
\begin{array}{cc}
     a^{2} & ab  \\
     ab & b^{2}
\end{array}
\right] .$$ 

\noindent  It is easy to check that

$$e_{1}e_{2}e_{1} = a^{2}e_{1}$$

\noindent and that 

$$e_{2}e_{1}e_{2} = a^{2}e_{2}.$$

\noindent Note also that  

$$e_{1}e_{2}=  \left[
\begin{array}{cc}
     a^{2} & ab  \\
     0 & 0
\end{array}
\right]$$ 

\noindent and 

$$e_{2}e_{1}=  \left[
\begin{array}{cc}
     a^{2} & 0  \\
     ab & 0
\end{array}
\right].$$

We define $$U_{i} = \delta e_{i}$$ \noindent for $i=1,2$ with $a^{2} = \delta^{-2}.$
Then we have , for $i = 1,2$

$$U_{i}^{2} = \delta U_{i}$$

$$U_{1}U_{2}U_{1} = U_{1}$$

$$U_{2}U_{1}U_{2} = U_{2}$$

\noindent and 

$$trace(U_{1})=trace(U_{2}) = \delta$$

\noindent while

$$trace(U_{1}U_{2}) = trace(U_{2}U_{1}) = 1.$$

\noindent We will use these results on the traces of these matrices in Section 6.
\vspace{3mm}

Now we return to the matrix parameters:
Since $a^{2} + b^{2} = 1$ this means that $\delta^{-2} + b^{2} = 1$ whence

$$b^{2} = 1-\delta^{-2}.$$

\noindent Therefore $b$ is real when $\delta^{2}$ is greater than or equal to $1$.
\vspace{3mm}

We are interested in the case where $\delta = -A^{2} - A^{-2}$ and {\em $A$ is a unit complex
number}.  Under these circumstances the braid group representation 

$$\rho(\sigma_{i}) = AU_{i} + A^{-1}1$$

\noindent will be unitary whenever $U_{i}$ is a real symmetric matrix. Thus we will obtain a
unitary representation of the three-strand braid group $B_{3}$ when $\delta^{2} \geq 1$.
Specifically, let $A=e^{i\theta}$. Then $\delta = -2cos(2\theta)$, so the condition 
$\delta^{2} \geq 1$ is equivalent to $cos^{2}(2\theta) \geq 1/4$. Thus we get the specific range of
angles $|\theta| \leq \pi/6$ and $|\theta - \pi| \leq \pi/6$ that gives unitary representations of
the three-strand braid group.  
\vspace{3mm}

\subsection{Pairs of Projectors and the Alexander Polynomial}
Just for the record we note a more general braid group representation that is available via our remarks about the structure of
two projectors.  Let $\{ W_{1},W_{2},...,W_{n-1},W_{n} \}$ be the standard basis of {\em column} vectors for a module of
dimension
$n$ over $k=C[A,A^{-1}]$ where $C$ denotes the complex numbers and 
$W_{k}$ is an $n$-tuple whose entries are zero in all places except the $k$-th place where the entry is one.
We shall refer to linear combinations of the $W_{k}$ as {\em vectors} over $k$.
Given any vector $v$ over $k$, let $|v>$ denote $v$ as a column vector, and let $<v| = v^{t}$ denote its transpose 
(just the transpose, as in our previous remarks), the
corresponding row vector.  Then $P(v) = |v><v|$ is a matrix such that $P^{2} = <v|v>P$, and $<v|v> = v^{t}v$ is equal to the
sum of the squares of the entries of $v.$ 
\vspace{3mm}

\noindent For $k = 1,2,..., n-1$ and $i^{2}=-1$, let
 
$$v_{k} = iAW_{k} - iA^{-1}W_{k+1}$$ \noindent and
 
$$U_{k} = |v_{k}><v_{k}|.$$ 

\noindent   Then, with $\delta = -A^{2} - A^{-2},$

$$U_{k}^{2} = \delta U_{k}$$
$$U_{k}U_{k \pm 1}U_{k} = U_{k}$$
$$U_{k}U_{l} = U_{l}U_{k}=0, |k-l|>1.$$

\noindent Thus these matrices give a special representation of the Temperley Lieb algebras $TL_{n}$ for each $n$.
Since the loop value is as given above, we can make correpsonding representations of the Artin Braid Groups $B_{n}$ by the
formulas

$$\rho(\sigma_{k}) = AI_{n} + A^{-1}U_{k},$$
$$\rho(\sigma_{k}^{-1}) = A^{-1}I_{n} + AU_{k}$$

\noindent where $I_{n}$ denotes the $n \times n$ identity matrix. It is not hard to verify that this representation of $B_{n}$
is equivalent to the classical Burau representation (See  \cite{KA91}) of the braid group. This
shows that there is a pathway from the algebra of projectors to the Alexander polynomial! We will treat this theme in a
separate paper.
\vspace{3mm}

\section{The Bracket Polynomial}
   
In this section we shall discuss the structure of the
the bracket state model for the Jones polynomial \cite{KA87}. In this way, we will explicitly construct the Jones polynomial by
using a state summation that is closely related to the braid group representation described in the last section. 
\vspace{3mm}

Before discussing the bracket polynomial we recall the basic theorem of Reidemeister \cite{Reid} about knot and link diagrams.
Reidemeister proved that the the three local moves on diagrams illustrated in Figure 6 capture combinatorially the notion of
ambient isotopy of links and knots in three-dimensional space. That is, if two diagrams represent knots or links that are 
isotopic in three-dimensional space, then the one diagram can be obtained from the other by a seqence of Reidemeister moves.
It is understood that a Reidemeister move is a local change on the diagram and that it is locally just as indicated by the
picture of the move. That is, a type one move adds or eliminates a loop in the underlying 4-regular graph of the knot
diagram. A type two move operates on a two sided region and a type three move operates on a three sided region. 
It is also understood that one can simplify a diagram by a homeomorphism of the plane. This could be called the type zero
move, but it is always available. The equivalence relation generated by the type two and type three moves is called {\em
regular isotopy}. The bracket polynomial is a regular isotopy invariant that can be normalized to produce an invariant of all
three Reidemeister moves.
\vspace{3mm}

{\tt    \setlength{\unitlength}{0.92pt}
\begin{picture}(357,406)
\thicklines   \put(215,374){\vector(-1,0){25}}
              \put(194,374){\vector(1,0){43}}
              \put(261,374){\line(1,0){86}}
              \put(97,360){\line(4,1){71}}
              \put(125,390){\line(-1,-1){29}}
              \put(52,372){\line(4,1){73}}
              \put(16,293){\framebox(18,19){1}}
              \put(110,308){\vector(1,0){43}}
              \put(134,308){\vector(-1,0){25}}
              \put(149,260){\line(3,2){57}}
              \put(163,341){\line(1,-1){43}}
              \put(72,308){\line(-1,1){21}}
              \put(95,288){\line(-5,4){15}}
              \put(96,316){\line(0,-1){27}}
              \put(38,280){\line(3,2){57}}
              \put(75,203){\line(-3,-5){26}}
              \put(66,178){\line(1,-1){15}}
              \put(47,198){\line(1,-1){12}}
              \put(20,69){\framebox(20,22){3}}
              \put(13,192){\framebox(24,23){2}}
              \put(10,372){\framebox(22,24){0}}
              \put(173,241){\line(0,-1){77}}
              \put(155,241){\line(0,-1){77}}
              \put(57,217){\line(-3,-5){10}}
              \put(71,242){\line(-1,-2){7}}
              \put(46,241){\line(3,-4){28}}
              \put(117,202){\vector(-1,0){25}}
              \put(94,202){\vector(1,0){43}}
              \put(55,130){\line(3,-4){28}}
              \put(80,131){\line(-1,-2){7}}
              \put(66,106){\line(-3,-5){10}}
              \put(83,93){\line(3,-4){28}}
              \put(108,94){\line(-1,-2){7}}
              \put(94,69){\line(-3,-5){10}}
              \put(57,50){\line(3,-4){28}}
              \put(85,53){\line(-1,-2){7}}
              \put(68,28){\line(-3,-5){10}}
              \put(222,131){\line(3,-4){28}}
              \put(247,132){\line(-1,-2){7}}
              \put(233,107){\line(-3,-5){10}}
              \put(199,90){\line(3,-4){28}}
              \put(224,91){\line(-1,-2){7}}
              \put(210,66){\line(-3,-5){10}}
              \put(226,54){\line(3,-4){28}}
              \put(251,53){\line(-1,-2){7}}
              \put(238,31){\line(-3,-5){10}}
              \put(56,89){\line(0,-1){38}}
              \put(108,94){\line(0,1){35}}
              \put(111,55){\line(0,-1){41}}
              \put(200,89){\line(0,1){43}}
              \put(200,48){\line(0,-1){36}}
              \put(250,94){\line(0,-1){43}}
              \put(136,77){\vector(1,0){43}}
              \put(159,77){\vector(-1,0){25}}
\end{picture}}

\begin{center}
{\bf Figure 6 - The Reidemeister Moves} 
\end{center}
\vspace{3mm}

The {\em bracket polynomial} , $<K> = <K>(A)$,  assigns to each unoriented
link diagram $K$ a Laurent polynomial in the variable $A$ such that
   
\begin{enumerate}
\item If $K$ and $K'$ are regularly isotopic links, then  $<K> = <K'>$.
  
\item If  $K \hspace{.1in} O$  denotes the disjoint union of $K$ with an extra
unknotted and unlinked component $O$, then
$$< K \hspace{.1in} O> = \delta<K>$$ where  $$\delta= -A^{2} - A^{-2}.$$
  
\item $<K>$ satisfies the following formula where in  Figure 7 the small diagrams represent
parts of larger diagrams that are identical except at the site indicated in the bracket.
In the text formula we have used the notations $S_{A}K$ and $S_{B}K$ to indicate the two smoothings
of a single crossing in the diagram $K$. That is, $K$,$S_{A}K$ and $S_{B}K$ differ at the site of
one crossing in the diagram $K$. The convention for these smoothings is indicated in Figure 7. 

$$<K> = A<S_{A}K> + A^{-1}<S_{B}K>$$

\end{enumerate}
\vspace{3mm}

{\tt    \setlength{\unitlength}{0.92pt}
\begin{picture}(427,204)
\thinlines    \put(275,8){\makebox(45,28){$+A$}}
              \put(118,7){\makebox(48,31){$=B$}}
              \put(118,87){\makebox(55,31){$=A$}}
              \put(275,90){\makebox(46,24){$+B$}}
\thicklines   \put(203,110){\line(1,0){39}}
              \put(242,112){\line(0,1){12}}
              \put(242,98){\line(0,-1){14}}
              \put(202,98){\line(1,0){40}}
              \put(202,84){\line(0,1){14}}
              \put(202,124){\line(0,-1){12}}
              \put(223,43){\line(1,0){13}}
              \put(222,3){\line(1,0){14}}
              \put(222,43){\line(0,-1){40}}
              \put(207,3){\line(-1,0){11}}
              \put(208,43){\line(0,-1){40}}
              \put(196,43){\line(1,0){12}}
              \put(42,2){\line(1,1){40}}
              \put(42,42){\line(1,-1){16}}
              \put(82,2){\line(-1,1){16}}
              \put(348,123){\line(1,0){12}}
              \put(360,123){\line(0,-1){40}}
              \put(359,83){\line(-1,0){11}}
              \put(374,123){\line(0,-1){40}}
              \put(374,83){\line(1,0){14}}
              \put(375,123){\line(1,0){13}}
              \put(357,43){\line(0,-1){12}}
              \put(357,3){\line(0,1){14}}
              \put(357,17){\line(1,0){40}}
              \put(397,17){\line(0,-1){14}}
              \put(397,31){\line(0,1){12}}
              \put(358,29){\line(1,0){39}}
              \put(349,44){\line(-1,-1){19}}
              \put(330,25){\line(1,-1){19}}
              \put(406,44){\line(1,-1){19}}
              \put(425,24){\line(-1,-1){19}}
              \put(190,44){\line(-1,-1){19}}
              \put(171,25){\line(1,-1){19}}
              \put(246,43){\line(1,-1){19}}
              \put(265,23){\line(-1,-1){19}}
              \put(38,44){\line(-1,-1){19}}
              \put(19,25){\line(1,-1){19}}
              \put(90,43){\line(1,-1){19}}
              \put(109,23){\line(-1,-1){19}}
              \put(343,123){\line(-1,-1){19}}
              \put(324,104){\line(1,-1){19}}
              \put(392,123){\line(1,-1){19}}
              \put(411,103){\line(-1,-1){19}}
              \put(196,125){\line(-1,-1){19}}
              \put(177,106){\line(1,-1){19}}
              \put(250,125){\line(1,-1){19}}
              \put(269,105){\line(-1,-1){19}}
              \put(47,122){\line(1,-1){40}}
              \put(47,82){\line(1,1){16}}
              \put(87,122){\line(-1,-1){16}}
              \put(41,123){\line(-1,-1){19}}
              \put(22,104){\line(1,-1){19}}
              \put(90,123){\line(1,-1){19}}
              \put(109,103){\line(-1,-1){19}}
              \put(304,161){\makebox(24,41){$=$}}
              \put(207,161){\makebox(50,41){$S_{B}$}}
              \put(337,201){\line(1,0){12}}
              \put(349,201){\line(0,-1){40}}
              \put(348,161){\line(-1,0){11}}
              \put(363,201){\line(0,-1){40}}
              \put(363,161){\line(1,0){14}}
              \put(364,201){\line(1,0){13}}
              \put(256,201){\line(1,-1){40}}
              \put(256,161){\line(1,1){16}}
              \put(296,201){\line(-1,-1){16}}
              \put(137,201){\line(0,-1){12}}
              \put(137,161){\line(0,1){14}}
              \put(137,175){\line(1,0){40}}
              \put(177,175){\line(0,-1){14}}
              \put(177,189){\line(0,1){12}}
              \put(138,187){\line(1,0){39}}
              \put(106,161){\makebox(24,41){$=$}}
              \put(1,160){\makebox(56,41){$S_{A}$}}
              \put(56,201){\line(1,-1){40}}
              \put(56,161){\line(1,1){16}}
              \put(96,201){\line(-1,-1){16}}
\end{picture}}

\begin{center}
{\bf Figure 7 - Bracket Expansion} 
\end{center}
\vspace{3mm}

\noindent This formula for expanding the bracket polynomial can be indicated symbolically in the same fashion that
we used in the previous section to indicate the representation of the Artin Braid Group to the Temperley Lieb algebra.
We will denote a crossing in the link diagram by the
letter chi, \mbox{\large $\chi$}. The
letter itself denotes a crossing where {\em the curved line in the letter chi is crossing over the straight segment in the
letter}. The barred letter denotes the switch of this crossing where {\em the curved line in the letter chi is undercrossing
the straight segment in the letter}. In the state model a crossing in a diagram for the
knot or link is expanded into two possible states by either smoothing (reconnecting) the crossing horizontally, \mbox{\large
$\asymp$}, or vertically $><$. Coefficients in this expansion correspond exactly to our representation of the braid group
so that any closed loop (without crossings) in the plane has value $\delta = -A^{2} - A^{-2}$ and the crossings expand
accrding to the formulas  
$$\mbox{\large $\chi$} = A \mbox{\large $\asymp$} + A^{-1} ><$$
$$\overline{\mbox{\large $\chi$}} = A^{-1} \mbox{\large $\asymp$} + A ><.$$
\noindent The verification that the bracket is invariant under the second Reidemeister move is then identical to
our proof in the previuos section that
$$\mbox{\large $\chi$}\overline{\mbox{\large $\chi$}} = \mbox{\large $\asymp$}.$$
\vspace{3mm}

{\tt    \setlength{\unitlength}{0.92pt}
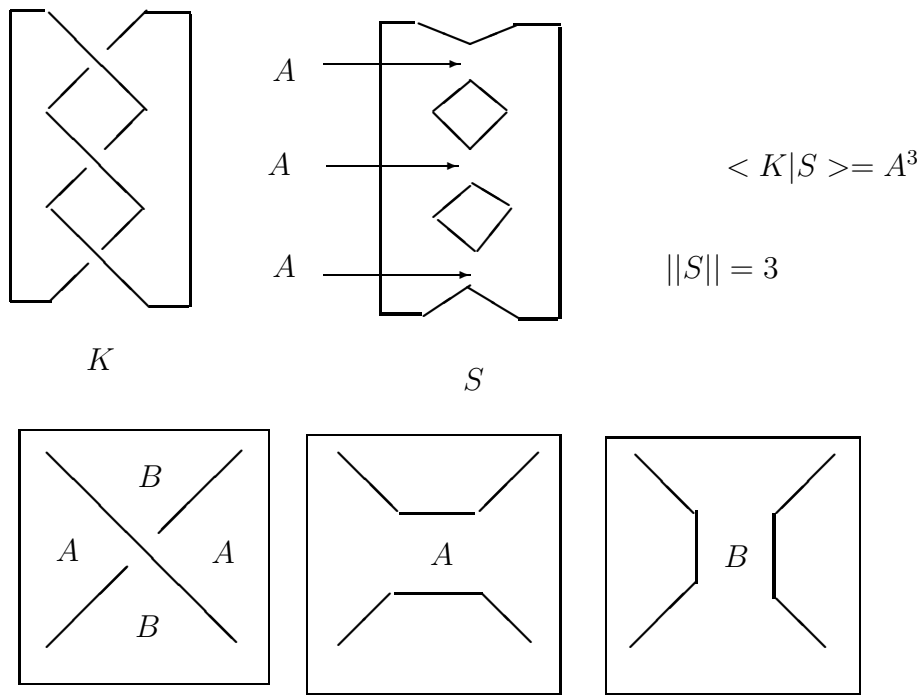
\begin{picture}(348,287)
\thicklines   \put(194,17){\line(3,-5){9}}
              \put(205,40){\makebox(86,47){$+ A^{-1}$}}
              \put(97,46){\makebox(44,40){$A$}}
              \put(67,67){\vector(1,0){21}}
              \put(202,199){\makebox(86,47){$+ A^{-1}$}}
              \put(94,205){\makebox(44,40){$A$}}
              \put(69,227){\vector(1,0){21}}
              \put(342,242){\line(0,1){43}}
              \put(174,46){\line(1,-5){4}}
              \put(196,43){\line(0,-1){19}}
              \put(345,47){\line(0,-1){41}}
              \put(323,47){\line(0,-1){41}}
              \put(180,17){\line(-1,-3){4}}
              \put(180,17){\line(1,0){15}}
              \put(180,25){\line(1,0){15}}
              \put(168,124){\line(3,-4){28}}
              \put(193,125){\line(-1,-2){7}}
              \put(179,100){\line(-3,-5){10}}
              \put(145,83){\line(3,-4){28}}
              \put(170,84){\line(-1,-2){7}}
              \put(156,59){\line(-3,-5){10}}
              \put(146,82){\line(0,1){43}}
              \put(146,41){\line(0,-1){36}}
              \put(196,87){\line(0,-1){43}}
              \put(317,127){\line(3,-4){28}}
              \put(342,128){\line(-1,-2){7}}
              \put(328,103){\line(-3,-5){10}}
              \put(294,86){\line(3,-4){28}}
              \put(319,87){\line(-1,-2){7}}
              \put(305,62){\line(-3,-5){10}}
              \put(295,85){\line(0,1){43}}
              \put(295,44){\line(0,-1){36}}
              \put(345,90){\line(0,-1){43}}
              \put(56,90){\line(0,-1){43}}
              \put(6,44){\line(0,-1){36}}
              \put(6,85){\line(0,1){43}}
              \put(44,27){\line(-3,-5){10}}
              \put(57,49){\line(-1,-2){7}}
              \put(32,50){\line(3,-4){28}}
              \put(16,62){\line(-3,-5){10}}
              \put(30,87){\line(-1,-2){7}}
              \put(5,86){\line(3,-4){28}}
              \put(39,103){\line(-3,-5){10}}
              \put(53,128){\line(-1,-2){7}}
              \put(28,127){\line(3,-4){28}}
              \put(317,241){\line(0,1){45}}
              \put(290,235){\line(0,1){51}}
              \put(192,277){\line(0,1){7}}
              \put(160,268){\line(1,2){8}}
              \put(141,267){\line(1,0){17}}
              \put(140,286){\line(0,-1){19}}
              \put(156,255){\line(3,-4){12}}
              \put(141,256){\line(1,0){15}}
              \put(140,236){\line(0,1){18}}
              \put(317,241){\line(3,-4){28}}
              \put(342,242){\line(-1,-2){7}}
              \put(328,217){\line(-3,-5){10}}
              \put(291,198){\line(3,-4){28}}
              \put(319,201){\line(-1,-2){7}}
              \put(302,176){\line(-3,-5){10}}
              \put(290,237){\line(0,-1){38}}
              \put(345,203){\line(0,-1){41}}
              \put(167,241){\line(3,-4){28}}
              \put(192,242){\line(-1,-2){7}}
              \put(178,217){\line(-3,-5){10}}
              \put(141,198){\line(3,-4){28}}
              \put(169,201){\line(-1,-2){7}}
              \put(152,176){\line(-3,-5){10}}
              \put(140,237){\line(0,-1){38}}
              \put(192,242){\line(0,1){35}}
              \put(195,203){\line(0,-1){41}}
              \put(58,207){\line(0,-1){41}}
              \put(55,246){\line(0,1){35}}
              \put(3,241){\line(0,-1){38}}
              \put(15,180){\line(-3,-5){10}}
              \put(32,205){\line(-1,-2){7}}
              \put(4,202){\line(3,-4){28}}
              \put(41,221){\line(-3,-5){10}}
              \put(55,246){\line(-1,-2){7}}
              \put(30,245){\line(3,-4){28}}
              \put(13,258){\line(-3,-5){10}}
              \put(27,283){\line(-1,-2){7}}
              \put(2,282){\line(3,-4){28}}
\end{picture}}

\begin{center}
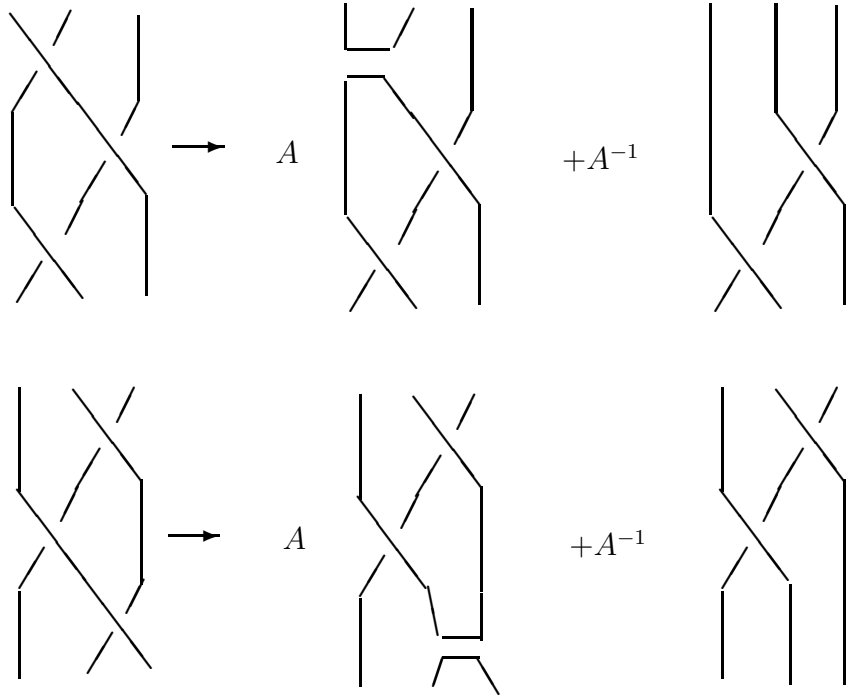
{\bf Figure 8 - Invariance of Bracket under Third Reidemeister Move} 
\end{center}
\vspace{3mm}

{\tt    \setlength{\unitlength}{0.92pt}
\begin{picture}(382,383)
\thicklines   \put(177,295){\line(6,5){15}}
              \put(194,281){\line(-6,5){15}}
\thinlines    \put(16,1){\makebox(242,89){$<K> = \Sigma_{S}<K|S>d^{||S||}$}}
              \put(249,249){\makebox(93,48){$||S|| = 3$}}
              \put(293,289){\makebox(88,53){$<K|S> = A^{3}$}}
              \put(98,259){\makebox(35,31){$A$}}
              \put(96,301){\makebox(35,31){$A$}}
              \put(98,340){\makebox(35,31){$A$}}
              \put(132,271){\vector(1,0){60}}
              \put(133,316){\vector(1,0){54}}
              \put(132,358){\vector(1,0){56}}
              \put(165,204){\makebox(57,48){$S$}}
              \put(20,214){\makebox(40,45){$K$}}
              \put(248,99){\framebox(104,105){}}
              \put(125,99){\framebox(104,106){}}
              \put(7,103){\framebox(102,104){}}
              \put(48,171){\makebox(25,35){$B$}}
              \put(47,109){\makebox(25,35){$B$}}
              \put(75,143){\makebox(31,28){$A$}}
              \put(11,144){\makebox(31,28){$A$}}
              \put(289,138){\makebox(25,35){$B$}}
              \put(165,142){\makebox(31,28){$A$}}
\thicklines   \put(285,174){\line(0,-1){30}}
              \put(258,118){\line(1,1){27}}
              \put(318,138){\line(1,-1){20}}
              \put(317,173){\line(0,-1){35}}
              \put(260,197){\line(1,-1){24}}
              \put(342,197){\line(-1,-1){24}}
              \put(138,119){\line(1,1){21}}
              \put(197,140){\line(-1,0){36}}
              \put(217,120){\line(-1,1){20}}
              \put(194,173){\line(-1,0){31}}
              \put(220,198){\line(-1,-1){24}}
              \put(138,198){\line(1,-1){24}}
              \put(98,199){\line(-1,-1){34}}
              \put(18,118){\line(1,1){33}}
              \put(18,198){\line(1,-1){78}}
              \put(208,300){\line(-5,3){15}}
              \put(195,281){\line(4,5){14}}
              \put(207,339){\line(-6,5){15}}
              \put(192,323){\line(1,1){15}}
              \put(177,338){\line(1,-1){15}}
              \put(192,351){\line(-6,-5){15}}
              \put(191,266){\line(5,-3){21}}
              \put(192,267){\line(-3,-2){19}}
              \put(192,366){\line(5,2){19}}
              \put(192,366){\line(-5,2){21}}
              \put(210,374){\line(1,0){19}}
              \put(229,373){\line(0,-1){120}}
              \put(212,253){\line(1,0){16}}
              \put(169,375){\line(-1,0){14}}
              \put(155,375){\line(0,-1){120}}
              \put(155,255){\line(1,0){17}}
              \put(3,260){\line(1,0){17}}
              \put(3,380){\line(0,-1){120}}
              \put(17,380){\line(-1,0){14}}
              \put(60,258){\line(1,0){16}}
              \put(77,378){\line(0,-1){120}}
              \put(58,379){\line(1,0){19}}
              \put(20,298){\line(1,-1){40}}
              \put(57,298){\line(-1,-1){16}}
              \put(19,260){\line(1,1){16}}
              \put(18,338){\line(1,-1){40}}
              \put(58,339){\line(-1,-1){16}}
              \put(18,299){\line(1,1){16}}
              \put(19,340){\line(1,1){16}}
              \put(59,380){\line(-1,-1){16}}
              \put(19,379){\line(1,-1){40}}
\end{picture}}

\begin{center}
{\bf Figure 9 - Bracket States} 
\end{center}
\vspace{3mm}
   
Knowing that the bracket is invariant under the second Reidemeister move allows us to verify directly that it is invariant
under the third Reidemeister move. This is illustrated in Figure 8.
In this Figure we show the two equivalent configurations in the third Reidemeister move vertically on the left, with
arrows point to the right of each configuration to an expansion via the bracket at one crossing.  The expansions give the same
bracket calculation due to invariance under the second Reidemeister move. Since the bracket
is invariant under the second and third Reidemeister moves, {\em property
1. is a direct consequence of properties 2. and 3.}. The second two properties
define the bracket on arbitrary link diagrams. 
\vspace{3mm}

  In fact we could have begun with the following more general definition: 
  Let $K$ be any unoriented link diaram. Define a 
  {\em state} of $K$  to be a choice of smoothings for all the crossings of $K.$  There are
  $2^{N}$ states of a diagram with $N$ crossings. A {\em smoothing} of a crossing
  is a local repacement of that crossing with two arcs that do not cross one another,
  as shown below. There are two choices for smoothing a given crossing. In illustrating
  a state it is convenient to label the smoothing with $A$ or $B$ to indicate the 
  crossing from which it was smoothed.  The $A$ or $B$ is called a {\em vertex
  weight} of the state.
  \vspace{3mm}

  Label each state with {\em vertex weights} $A$ or $B$ as illustrated in Figure 9. Here $A$ and
  $B$ are commuting polynomial variables. Define two evaluations related to the state:
  The first evaluation is the product of the vertex weights, denoted  $$[K|S].$$
  The second evaluation is the number of loops (Jordan curves) 
  in the state $S$, denoted  $$||S||.$$
  
  Define the {\em state summation}, $[K]$, by the formula 
  $$[K]=\sum_{S} [K|S]\delta^{||S||-1}.$$
  It follows from this definition,that $[K]$ satisfies the formulas
  
$$[\mbox{\large $\chi$}] = A [\mbox{\large $\asymp$}] + B [><]$$
  $$[O \hspace{.1in} K] = \delta[K],$$ and 
  $$[O] =1.$$
  
  The demand that
  $[K]$ be invariant under the {\em second} Reidemeister move leads to the conditions
  $B=A^{-1}$ and $\delta=-A^{2}+A^{-2}$. This specialization is easily seen to be 
  invariant under the third Reidemeister move. Calling this specialization the {\em 
  topological} bracket, and denoting it (as above) by $<K>$ one finds the following 
  behaviour
  under the first Reidemeister move 
  $$<\mbox{\large $\gamma$}> = -A^{3} <\smile>$$ and 
  $$<\overline{\mbox{\large $\gamma$}}> = -A^{-3} <\smile>$$

\noindent where \mbox{\large $\gamma$}  denotes a curl of positive type as indicated in Figure 10, and 
$\overline{\mbox{\large $\gamma$}}$ indicates a curl of negative type as also seen in this Figure.

  The topological bracket is invariant under regular isotopy and can be  
  normalized to an invariant of ambient isotopy by the definition  
  $$f_{K}(A) = (-A^{3})^{-w(K)}<K>(A)$$ where 
  w(K) is the sum of the crossing signs  of the oriented
  link K. w(K) is called the writhe of K.  The convention for crossing signs is shown
  in  Figure 10.
  \vspace{3mm}

{\tt    \setlength{\unitlength}{0.92pt}
\begin{picture}(282,126)
\thicklines   \put(226,17){\vector(1,-1){14}}
              \put(200,42){\line(1,-1){17}}
              \put(200,2){\vector(1,1){40}}
              \put(64,26){\vector(1,1){16}}
              \put(40,3){\line(1,1){17}}
              \put(40,43){\vector(1,-1){40}}
              \put(242,102){\line(1,0){14}}
              \put(89,103){\line(1,0){13}}
              \put(95,109){\line(0,-1){12}}
              \put(91,27){\line(0,-1){12}}
              \put(85,21){\line(1,0){13}}
              \put(240,21){\line(1,0){14}}
              \put(201,83){\line(1,0){39}}
              \put(241,123){\line(1,0){40}}
              \put(160,123){\line(1,0){41}}
              \put(240,83){\line(-1,1){15}}
              \put(200,123){\line(1,-1){16}}
              \put(240,123){\line(-1,-1){39}}
              \put(81,122){\line(1,0){40}}
              \put(1,123){\line(1,0){39}}
              \put(41,82){\line(1,0){38}}
              \put(41,83){\line(1,1){14}}
              \put(80,123){\line(-1,-1){15}}
              \put(40,122){\line(1,-1){38}}
\end{picture}}

\begin{center}
{\bf Figure 10 - Crossing Signs} 
\end{center}
\vspace{3mm}

  By a change of variables one obtains the original
  Jones polynomial, $V_{K}(t)$ \cite{JO86} from the normalized bracket:

  $$V_{K}(t) = f_{K}(t^{-1/4}).$$

{\tt    \setlength{\unitlength}{0.92pt}
\begin{picture}(196,168)
\thicklines   \put(130,1){\makebox(55,38){$T^{*}$}}
              \put(19,6){\makebox(44,33){$T$}}
              \put(193,45){\line(-1,0){14}}
              \put(193,165){\line(0,-1){121}}
              \put(180,165){\line(1,0){14}}
              \put(125,46){\line(1,0){15}}
              \put(125,165){\line(0,-1){119}}
              \put(140,165){\line(-1,0){15}}
              \put(180,45){\line(-1,1){16}}
              \put(180,85){\line(-1,1){15}}
              \put(180,124){\line(-1,1){16}}
              \put(140,85){\line(1,-1){16}}
              \put(141,125){\line(1,-1){16}}
              \put(140,165){\line(1,-1){18}}
              \put(179,85){\line(-1,-1){39}}
              \put(180,125){\line(-1,-1){40}}
              \put(180,165){\line(-1,-1){39}}
              \put(5,45){\line(1,0){17}}
              \put(3,165){\line(0,-1){120}}
              \put(19,165){\line(-1,0){14}}
              \put(62,43){\line(1,0){16}}
              \put(79,163){\line(0,-1){120}}
              \put(60,164){\line(1,0){19}}
              \put(22,83){\line(1,-1){40}}
              \put(59,83){\line(-1,-1){16}}
              \put(21,45){\line(1,1){16}}
              \put(20,123){\line(1,-1){40}}
              \put(60,124){\line(-1,-1){16}}
              \put(20,84){\line(1,1){16}}
              \put(21,125){\line(1,1){16}}
              \put(61,165){\line(-1,-1){16}}
              \put(21,164){\line(1,-1){40}}
\end{picture}}

\begin{center}
{\bf Figure 11 - Trefoil and Mirror Image} 
\end{center}
\vspace{3mm}

The bracket model for
the Jones polynomial is quite useful both theoretically and in terms of practical
computations. One of the neatest applications is to simply compute $f_{K}(A)$ for the
trefoil knot $T$ and determine that $f_{K}(A)$ is not equal to $f_{K}(A^{-1}).$  This shows that
the trefoil is not ambient isotopic to its mirror image (See Figure 11), a fact that is quite tricky
to prove by classical methods.
\vspace{3mm}

\noindent {\bf Remark.}
The relationship of the Temperley Lieb algebra  with the bracket polynomial comes through the basic
bracket identity. This identity, interpreted in the context of the diagrammtic Temperley Lieb algebra becomes a
representation $\rho$ of the Artin braid group $B_{n}$ on $n$ strands to the Temperley Lieb algebra
$TL_{n}$ defined by the formulas

$$\rho(\sigma_{i}) = AU_{i} + A^{-1}1$$

$$\rho(\sigma_{i}^{-1}) = A^{-1}U_{i} + A1.$$

\noindent Here $sigma_{i}$ denotes the braid generator that twists strands $i$ and $i+1$. For this representation of
the Temperley Lieb algebra, the loop value $\delta$ is
$-A^{2}-A^{-2}$ and the ring $k$ is $Z[A,A^{-1}]$, the ring of Laurent polynomials in $A$ with
integer coefficients. \vspace{3mm}

\noindent {\bf Remark.}
There are hints of quantum mechanical interpretations in the combinatorics of this state sum model
for the Jones polynomial. The expansion formula for the bracket polynomial

$$<K> = A<S_{A}K> + A^{-1}<S_{B}K>$$

\noindent suggests that the diagram $K$ should be thought of as a superposition of the diagrams
$S_{A}K$ and $S_{B}K$. That is, we can think of a knot diagram with respect to a given crossing as
the superposition of the diagrams obtained by smoothing that crossing.  Then, with respect to all
the crossings, one can think of the diagram as a superposition of the states obtained by smoothing
each crossing in one of its two possible ways. This is a superposition view of the bracket state
sum as a whole.

$$<K>=\sum_{S} <K|S>\delta^{||S||-1}$$

\noindent In this sense the bracket polynomial evaluation is directly analogous to an amplitude in
quantum mechanics. We shall make this analogy more precise in the sections to follow. However, the
topological information is contained in this amplitude as whole, and not in any specific state
evaluation. Thus the topological model ignores the standard measurement situation in quantum
mechanics where one gets at best information about one state at a time when a measurement is taken.
This means that a quantum computational model of the bracket polynomial will be essentially
probabilistic, only giving partial information at each measurement.
\vspace{3mm}

As a result of this discussion, it is natural to ask to what extent one can extract partial
topological information from an incomplete summation over the states of the bracket polynomial.
It is not clear at this stage what this answer is to this question. It may require a new
exploration of the properties of the state sum and its corresponding polynomial.
\vspace{3mm}

\section{Knot Amplitudes}
 
At the end of the first section we said: the connection of quantum
mechanics with  topology is an amplification of Dirac notation. In this section we begin the process of amplification!
\vspace{3mm}

{\tt    \setlength{\unitlength}{0.92pt}
\begin{picture}(154,111)
\thicklines   \put(77,60){\circle{40}}
\thinlines    \put(1,78){\makebox(33,32){$t$}}
              \put(127,5){\makebox(26,26){$x$}}
              \put(25,20){\vector(1,0){92}}
              \put(38,1){\vector(0,1){99}}
\end{picture}}

\begin{center}
{\bf Figure 12 - Circle in Spacetime} 
\end{center}
\vspace{3mm}

Consider first a circle in a spacetime plane with time represented vertically and
space horizontally. The circle represents a vacuum to vacuum process that
includes the creation of
two "particles", and their subsequent annihilation. See Figures 12 and 13.
\vspace{20mm}

 {\tt    \setlength{\unitlength}{0.92pt}
\begin{picture}(104,132)
\thicklines   \put(3,1){\vector(0,1){130}}
\thinlines    \put(63,44){\oval(80,78)[b]}
              \put(63,84){\oval(80,78)[t]}
\end{picture}}
 
\begin{center}
{\bf Figure 13 - Creation and Annihilation} 
\end{center}
\vspace{3mm}

In accord with our previous description, we could divide the circle into these
two parts (creation(a)  and annihilation (b)) and consider the amplitude   $<b|a>.$
Since the diagram for the creation of the two particles ends in two separate
points, it is natural to take a vector space of the form  $V \otimes V$  as the target for
the bra and as the domain of the ket.
\vspace{3mm}

We imagine at least one particle property being catalogued by each dimension of
$V.$  For example, a basis of  $V$  could enumerate the spins of the created
particles.  If $\{e_{a} \}$ is a basis for $V$ then $\{e_{a} \otimes e_{b} \}$ forms a basis for
$V \otimes V.$ The elements of this new basis constitute all possible combinations of the
particle properties. Since such combinations are multiplicative, the tensor product is the
appropriate construction.
\vspace{3mm}

In this language the creation ket is a map  $cup$,

$$cup = |a> : C \longrightarrow V \otimes V,$$

\noindent and the annihilation bra is a mapping  $cap$,

$$cap= <b| : V \otimes V \longrightarrow C.$$

The first hint of topology comes when we realise that it is possible to draw a
much more complicated simple closed curve in the plane that is nevertheless
decomposed with respect to the vertical direction into many cups and caps.  In
fact, any simple (no self-intersections) differentiable curve can be rigidly rotated until
it is in general position with respect to the vertical.  It will then be seen to
be decomposed into these minima and maxima.   Our prescriptions for amplitudes
suggest that we regard any such curve as an amplitude via its description as a
mapping from   $C$  to $C.$
\vspace{3mm}

Each simple closed curve gives rise to an amplitude,  but any simple closed curve
in the plane is isotopic to a circle, by the Jordan Curve Theorem.  If these are
topological amplitudes,  then they should all be equal to the original amplitude
for the circle.   Thus the question:  What condition on creation and annihilation
will insure topological amplitudes?  The answer derives from the fact that all
isotopies of the simple closed curves are generated by the cancellation of
adjacent  maxima and minima as illustrated below.
\vspace{3mm}

{\tt    \setlength{\unitlength}{0.92pt}
\begin{picture}(206,181)
\thinlines    \put(203,173){\line(0,-1){172}}
              \put(169,94){\vector(-1,0){45}}
              \put(138,94){\vector(1,0){48}}
              \put(133,107){\line(-1,-1){1}}
              \put(109,103){\line(0,1){77}}
              \put(3,78){\line(0,-1){76}}
              \put(83,101){\oval(54,86)[b]}
              \put(30,76){\oval(54,86)[t]}
\end{picture}}
  
\begin{center}
{\bf Figure 14 - Cancellation of Maxima and Minima} 
\end{center}
\vspace{3mm}

In composing mappings it is necessary to use the identifications 
$(V \otimes V) \otimes V = V \otimes (V \otimes V)$ 
and $V \otimes k=k \otimes V = V.$ Thus in the illustration above, the composition on the left
is given by

$$V = V \otimes k  - 1 \otimes cup \rightarrow  V \otimes (V \otimes V)$$

$$ = (V \otimes V) \otimes V  - cap \otimes 1 \rightarrow  k \otimes V = V.$$

\noindent This composition must equal the identity map on $V$ (denoted $1$ here) for the
amplitudes to have a proper image of the topological cancellation.
This condition is said very simply by taking a matrix representation for the
corresponding operators.
\vspace{3mm}

Specifically,  let   $\{e_{1}, e_{2}, ..., e_{n} \}$  be a basis for $V.$ 
Let $e_{ab} = e_{a} \otimes  e_{b}$  
denote the elements of the tensor basis for $V \otimes V.$  Then there are matrices  $M_{ab}$  
and  $M^{ab}$  such that

$$cup(1)  =  \Sigma M^{ab} e_{ab}$$

\noindent   with the summation taken over all values of $a$ and $b$ 
from $1$ to $n.$  Similarly,  $cap$  is described by 

$$cap(e_{ab}) =  M_{ab.}$$  

\noindent Thus the
amplitude for the circle is 

$$cap[cup(1)]  =  cap \Sigma M^{ab}e_{ab} = \Sigma M^{ab}M_{ab.}$$

\noindent In
general, the value of the amplitude on a simple closed curve is obtained by translating it into an
``abstract tensor expression"  in the $M_{ab}$ and $M^{ab}$, and then summing over these
products for all cases of repeated indices.
\vspace{3mm}

Returning to the topological conditions we see that they are just that the
matrices  $(M_{ab})$  and  $(M^{ab})$  are inverses in the sense that   
$\Sigma M_{ai}M^{ib}  = \delta_{a}^{b}$  and
$\Sigma M^{ai}M_{ib}  = \delta^{a}_{b}$ where $\delta_{a}^{b}$ denotes the (identity matrix)
Kronecker delta that is equal to one when its two indices are equal to one another and zero
otherwise.
\vspace{3mm}

 {\tt    \setlength{\unitlength}{0.92pt}
\begin{picture}(247,219)
\thinlines    \put(216,172){\makebox(28,35){$b$}}
              \put(217,3){\makebox(29,31){$a$}}
              \put(119,183){\makebox(28,35){$b$}}
              \put(67,88){\makebox(23,30){$i$}}
              \put(17,1){\makebox(29,31){$a$}}
              \put(206,188){\circle*{10}}
              \put(6,18){\circle*{10}}
              \put(60,103){\circle*{10}}
              \put(113,196){\circle*{10}}
              \put(-24,206){\circle*{0}}
              \put(206,17){\circle*{10}}
              \put(206,189){\line(0,-1){172}}
              \put(172,110){\vector(-1,0){45}}
              \put(141,110){\vector(1,0){48}}
              \put(136,123){\line(-1,-1){1}}
              \put(112,119){\line(0,1){77}}
              \put(6,94){\line(0,-1){76}}
              \put(86,117){\oval(54,86)[b]}
              \put(33,92){\oval(54,86)[t]}
\end{picture}}

\begin{center}
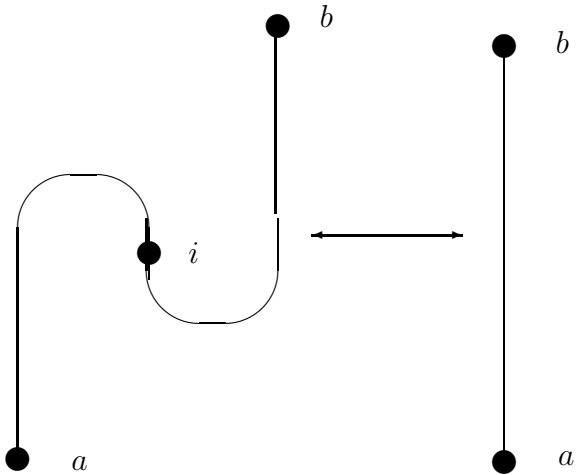
{\bf Figure 15 - Algebraic Cancellation of Maxima and Minima} 
\end{center}
\vspace{3mm}

In Figure 15, we show the diagrammatic representative of  the
equation $\Sigma M_{ai}M^{ib}  = \delta_{a}^{b}.$
\vspace{3mm}

In the simplest case  $cup$ and $cap$  are represented by  $2 \times 2$ matrices.  The
topological condition implies that these matrices are inverses of each other.
Thus the problem of the existence of topological amplitudes is very easily solved
for simple closed curves in the plane.
\vspace{3mm}

Now we go to knots and links.   Any knot or link can be represented by a picture
that is configured with respect to a vertical direction in the plane.  The
picture will decompose into minima (creations)  maxima (annihilations)  and
crossings of the two types shown below.  (Here I consider knots and links that
are unoriented.  They do not have an intrinsic preferred direction of travel.)
See Figure 16. In Figure 16 we have indicated the crossings as  
mappings of  $V \otimes V$  to itself ,  called   $R$  and  $R^{-1}$  respectively.  These
mappings represent the transitions corresponding to these elementary configurations.
\vspace{3mm}

{\tt    \setlength{\unitlength}{0.92pt}
\begin{picture}(348,219)
\thinlines    \put(294,58){\makebox(53,46){$R^{-1}$}}
              \put(296,138){\makebox(43,47){$R$}}
\thicklines   \put(195,50){\line(0,1){18}}
              \put(156,51){\line(1,0){39}}
              \put(154,66){\line(0,-1){15}}
              \put(118,51){\line(0,1){15}}
              \put(76,52){\line(1,0){42}}
              \put(75,66){\line(0,-1){14}}
              \put(156,162){\line(0,-1){18}}
              \put(118,163){\line(1,0){38}}
              \put(118,148){\line(0,1){15}}
              \put(278,64){\circle*{12}}
              \put(239,103){\circle*{12}}
              \put(238,65){\circle*{12}}
              \put(276,104){\circle*{12}}
              \put(276,145){\circle*{12}}
              \put(238,184){\circle*{12}}
              \put(239,146){\circle*{12}}
              \put(276,184){\circle*{12}}
\thinlines    \put(1,183){\makebox(32,35){$t$}}
              \put(228,1){\makebox(60,47){$x$}}
              \put(25,27){\vector(1,0){193}}
              \put(36,11){\vector(0,1){194}}
\thicklines   \put(157,106){\circle*{12}}
              \put(197,67){\circle*{12}}
              \put(77,68){\circle*{12}}
              \put(75,165){\circle*{12}}
              \put(119,67){\circle*{12}}
              \put(155,67){\circle*{12}}
              \put(115,107){\circle*{12}}
              \put(118,145){\circle*{12}}
              \put(156,145){\circle*{12}}
              \put(196,166){\circle*{12}}
              \put(196,187){\line(-1,0){40}}
              \put(196,147){\line(0,1){40}}
              \put(75,187){\line(1,0){41}}
              \put(75,146){\line(0,1){41}}
              \put(196,107){\line(0,-1){42}}
              \put(76,105){\line(0,-1){38}}
              \put(116,187){\line(1,0){40}}
              \put(276,105){\line(-1,-1){39}}
              \put(237,105){\line(1,-1){16}}
              \put(277,65){\line(-1,1){16}}
              \put(237,146){\line(1,1){16}}
              \put(277,186){\line(-1,-1){16}}
              \put(237,185){\line(1,-1){40}}
              \put(155,67){\line(-1,1){16}}
              \put(116,104){\line(1,-1){16}}
              \put(156,106){\line(-1,-1){39}}
              \put(156,146){\line(1,-1){40}}
              \put(196,147){\line(-1,-1){16}}
              \put(156,107){\line(1,1){16}}
              \put(76,105){\line(1,1){16}}
              \put(116,145){\line(-1,-1){16}}
              \put(76,144){\line(1,-1){40}}
\end{picture}}

\begin{center}
{\bf Figure 16 - Morse Knot Decomposition} 
\end{center}
\vspace{3mm}

That  $R$  and $R^{-1}$  really must be inverses follows from the isotopy shown in Figure 17
(This is the second Reidemeister move.)
\vspace{3mm}

{\tt    \setlength{\unitlength}{0.92pt}
\begin{picture}(272,166)
\thicklines   \put(226,5){\makebox(45,37){$d$}}
              \put(181,8){\makebox(44,33){$c$}}
              \put(228,134){\makebox(39,28){$b$}}
              \put(178,134){\makebox(44,26){$a$}}
              \put(89,69){\makebox(33,40){$j$}}
              \put(1,74){\makebox(34,34){$i$}}
              \put(66,1){\makebox(38,38){$d$}}
              \put(21,1){\makebox(41,38){$c$}}
              \put(63,133){\makebox(36,32){$b$}}
              \put(25,133){\makebox(35,30){$a$}}
              \put(153,89){\vector(-1,0){18}}
              \put(135,89){\vector(1,0){40}}
              \put(202,128){\circle*{12}}
              \put(201,48){\circle*{12}}
              \put(241,128){\line(0,-1){80}}
              \put(201,129){\line(0,-1){80}}
              \put(83,45){\circle*{12}}
              \put(241,48){\circle*{12}}
              \put(43,46){\circle*{12}}
              \put(241,128){\circle*{12}}
              \put(81,87){\circle*{12}}
              \put(43,126){\circle*{12}}
              \put(44,88){\circle*{12}}
              \put(81,126){\circle*{12}}
              \put(81,86){\line(-1,-1){39}}
              \put(42,86){\line(1,-1){16}}
              \put(82,46){\line(-1,1){16}}
              \put(42,88){\line(1,1){16}}
              \put(82,128){\line(-1,-1){16}}
              \put(42,127){\line(1,-1){40}}
\end{picture}}

\begin{center}
{\bf Figure 17 - Braiding Cancellation} 
\end{center}
\vspace{3mm}

We now have the vocabulary of  $cup$,$cap$,  $R$  and $R^{-1}$.   Any knot or link can be
written as a composition of these fragments, and consequently a choice of such
mappings determines an amplitude for knots and links.  In order for such an
amplitude to be topological    we want it to be invariant under the list of local
moves on the diagrams shown in Figure 18.  These moves are an augmented list of
the Reidemeister moves, adjusted to take care of the fact that the diagrams are
arranged with respect to a given direction in the plane. The equivalence relation
generated by these moves is called regular isotopy.  It is one move short of the
relation known as ambient isotopy.  The missing move is the first Reidemeister
move shown in Figure 6.
\vspace{3mm}

In the first Reidemeister move, a curl in the diagram is created or destroyed. 
Ambient isotopy (generated by all the Reidemeister moves) corresponds to the full
topology of knots and links embedded in three dimensional space. Two link
diagrams are ambient isotopic via the Reidemeister moves if and only if there is
a continuous family of embeddings in three dimensions leading from one link to
the other.  The moves give us a combinatorial reformulation of the spatial
topology of knots and links.
\vspace{3mm}

{\tt    \setlength{\unitlength}{0.92pt}
\begin{picture}(260,421)
\thicklines   \put(2,43){\framebox(40,41){$4$}}
              \put(1,163){\framebox(42,40){$3$}}
              \put(2,257){\framebox(43,42){$2$}}
              \put(1,363){\framebox(42,41){$0$}}
              \put(131,24){\vector(1,0){43}}
              \put(154,24){\vector(-1,0){25}}
              \put(248,37){\line(-1,1){36}}
              \put(237,25){\line(1,1){11}}
              \put(212,2){\line(1,1){14}}
              \put(232,41){\line(0,-1){40}}
              \put(192,42){\line(1,0){40}}
              \put(192,2){\line(0,1){40}}
              \put(53,38){\line(1,1){36}}
              \put(68,22){\line(-1,1){15}}
              \put(89,3){\line(-1,1){12}}
              \put(111,43){\line(0,-1){40}}
              \put(72,43){\line(1,0){39}}
              \put(72,3){\line(0,1){40}}
              \put(109,374){\vector(1,0){43}}
              \put(132,374){\vector(-1,0){25}}
              \put(162,420){\line(0,-1){95}}
              \put(101,351){\line(0,1){69}}
              \put(78,350){\line(1,0){22}}
              \put(78,380){\line(0,-1){30}}
              \put(53,380){\line(1,0){25}}
              \put(53,324){\line(0,1){56}}
              \put(78,263){\line(-3,-5){26}}
              \put(69,238){\line(1,-1){15}}
              \put(50,258){\line(1,-1){12}}
              \put(176,301){\line(0,-1){77}}
              \put(158,301){\line(0,-1){77}}
              \put(60,277){\line(-3,-5){10}}
              \put(74,302){\line(-1,-2){7}}
              \put(49,301){\line(3,-4){28}}
              \put(120,262){\vector(-1,0){25}}
              \put(97,262){\vector(1,0){43}}
              \put(58,203){\line(3,-4){28}}
              \put(83,204){\line(-1,-2){7}}
              \put(69,179){\line(-3,-5){10}}
              \put(86,166){\line(3,-4){28}}
              \put(111,167){\line(-1,-2){7}}
              \put(97,142){\line(-3,-5){10}}
              \put(60,123){\line(3,-4){28}}
              \put(88,126){\line(-1,-2){7}}
              \put(71,101){\line(-3,-5){10}}
              \put(225,204){\line(3,-4){28}}
              \put(250,205){\line(-1,-2){7}}
              \put(236,180){\line(-3,-5){10}}
              \put(202,163){\line(3,-4){28}}
              \put(227,164){\line(-1,-2){7}}
              \put(213,139){\line(-3,-5){10}}
              \put(229,127){\line(3,-4){28}}
              \put(254,126){\line(-1,-2){7}}
              \put(241,104){\line(-3,-5){10}}
              \put(59,162){\line(0,-1){38}}
              \put(111,167){\line(0,1){35}}
              \put(114,128){\line(0,-1){41}}
              \put(203,162){\line(0,1){43}}
              \put(203,121){\line(0,-1){36}}
              \put(253,167){\line(0,-1){43}}
              \put(139,150){\vector(1,0){43}}
              \put(162,150){\vector(-1,0){25}}
\end{picture}}

\begin{center}
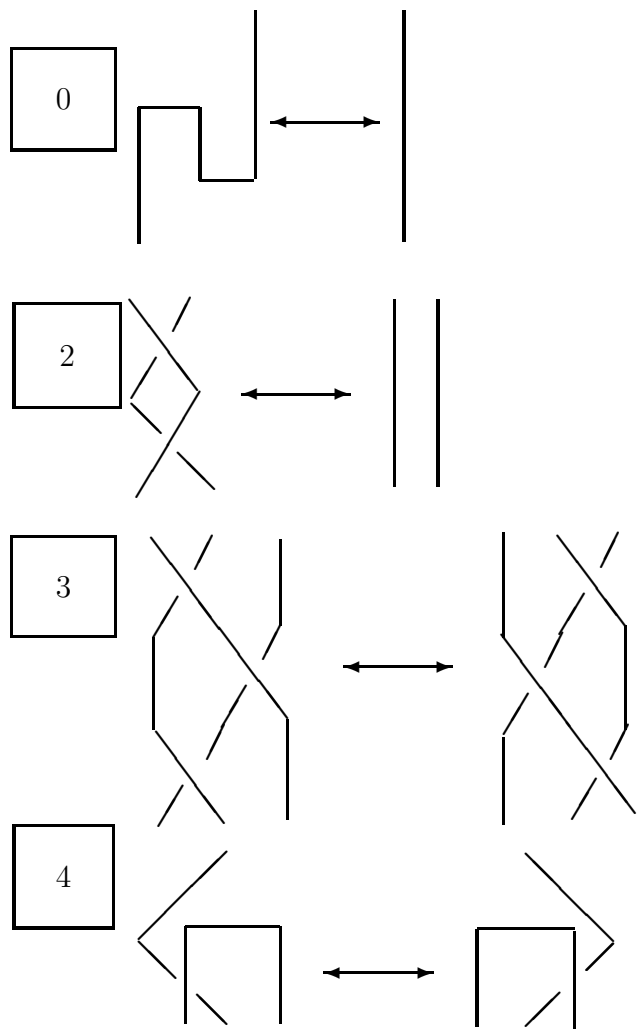
{\bf Figure 18- Moves for Regular Isotopy of Morse Diagrams} 
\end{center}
\vspace{3mm}

By ignoring the first Reidemeister move, we allow the possibility that these
diagrams can model framed links, that is links with a normal vector field
or,equivalently, embeddings of curves that are thickened into bands. It turns out
to be fruitful to study invariants of regular isotopy. In fact, one can usually
normalise an invariant of regular isotopy to obtain an invariant of ambient
isotopy. We shall see an example of this phenomenon with the bracket polynomial
in a few paragraphs.
\vspace{3mm}

As the reader can see, we have already discussed the algebraic meaning of moves 
0. and 2.  The other moves translate into very interesting algebra.  Move 3.,
when translated into algebra, is the  famous Yang-Baxter equation.    The
Yang-Baxter equation occurred for the first time in problems related to  exactly
solved models in statistical mechanics (See \cite{KA96}.). All the moves taken
together are directly related to the axioms for a quasi-triangular Hopf algebra
(aka quantum group).  We shall not go into this connection here.
\vspace{3mm}

There is an intimate connection between knot invariants and the structure of
generalised amplitudes, as we have described them in terms of vector space
mappings associated with link diagrams.  This strategy for the construction of 
invariants is directly motivated by the concept of an amplitude in quantum
mechanics.  It turns out that the invariants that can actually be produced by
this means (that is by assigning finite dimensional matrices to the caps, cups
and crossings) are incredibly rich.  They encompass, at present,  all of the
known invariants of polynomial type  (Alexander polynomial, Jones polynomial and
their generalisations.).
\vspace{3mm}

It is now possible to indicate the construction of the Jones polynomial via the
bracket polynomial as an amplitude, by specifying its matrices. The cups and the
caps are defined by  $(M_{ab}) = (M^{ab}) = M$  where $M$ is the $2 \times 2$ matrix (with
$ii=-1$).\vspace{3mm}

$$M =  \left[
\begin{array}{cc}
     0 & iA  \\
     -iA^{-1} & 0
\end{array}
\right] $$

\noindent Note that  $MM = I$  where  $I$ is the identity matrix.   Note also that the
amplitude for the circle is

$$\Sigma M_{ab}M^{ab} = \Sigma M_{ab} M_{ab} = \Sigma M_{ab}^{2}$$

$$=  (iA)^{2}  +  (-iA^{-1})^{2}  =  - A^{2}  - A^{-2.}$$

\noindent The matrix   $R$   is then defined by the equation

$$R^{ab}_{cd}   =  AM^{ab}M_{cd}  +  A^{-1} \delta^{a}_{c} \delta^{b}_{d},$$
ªª
 
\noindent Since, diagrammatically,  we identify  $R$ with a (right handed) crossing, this
equation can be written diagrammatically as the generating identity for the bracket polynomial:
\vspace{3mm}

$$\mbox{\large $\chi$} = A \mbox{\large $\asymp$} + A^{-1} ><$$

{\tt    \setlength{\unitlength}{0.92pt}
\begin{picture}(335,308)
\thicklines   \put(83,1){\makebox(240,78){$\eta^{b}_{a} = M_{ai}M^{bi}$}}
              \put(296,176){\makebox(38,44){$i$}}
              \put(203,138){\makebox(29,36){$b$}}
              \put(202,216){\makebox(33,32){$a$}}
              \put(243,157){\circle*{12}}
              \put(242,237){\circle*{12}}
              \put(284,198){\circle*{12}}
              \put(243,118){\line(0,1){39}}
              \put(282,118){\line(-1,0){39}}
              \put(283,278){\line(0,-1){160}}
              \put(243,280){\line(1,0){40}}
              \put(241,238){\line(0,1){42}}
              \put(3,91){\line(0,1){31}}
              \put(163,91){\line(-1,0){161}}
              \put(163,305){\line(0,-1){212}}
              \put(3,305){\line(1,0){159}}
              \put(3,280){\line(0,1){25}}
              \put(42,106){\line(0,1){14}}
              \put(144,106){\line(-1,0){102}}
              \put(144,292){\line(0,-1){186}}
              \put(43,293){\line(1,0){101}}
              \put(43,279){\line(0,1){14}}
              \put(121,120){\line(-1,0){40}}
              \put(122,279){\line(0,-1){159}}
              \put(83,279){\line(1,0){39}}
              \put(3,159){\line(0,-1){39}}
              \put(82,198){\line(0,-1){38}}
              \put(82,280){\line(0,-1){42}}
              \put(3,238){\line(0,-1){38}}
              \put(83,120){\line(-1,1){16}}
              \put(82,198){\line(-1,1){16}}
              \put(44,158){\line(1,-1){16}}
              \put(42,239){\line(1,-1){18}}
              \put(81,159){\line(-1,-1){39}}
              \put(82,239){\line(-1,-1){39}}
              \put(4,198){\line(1,-1){40}}
              \put(42,198){\line(-1,-1){16}}
              \put(4,160){\line(1,1){16}}
              \put(3,239){\line(1,1){16}}
              \put(43,279){\line(-1,-1){16}}
              \put(3,278){\line(1,-1){40}}
\end{picture}}

\begin{center}
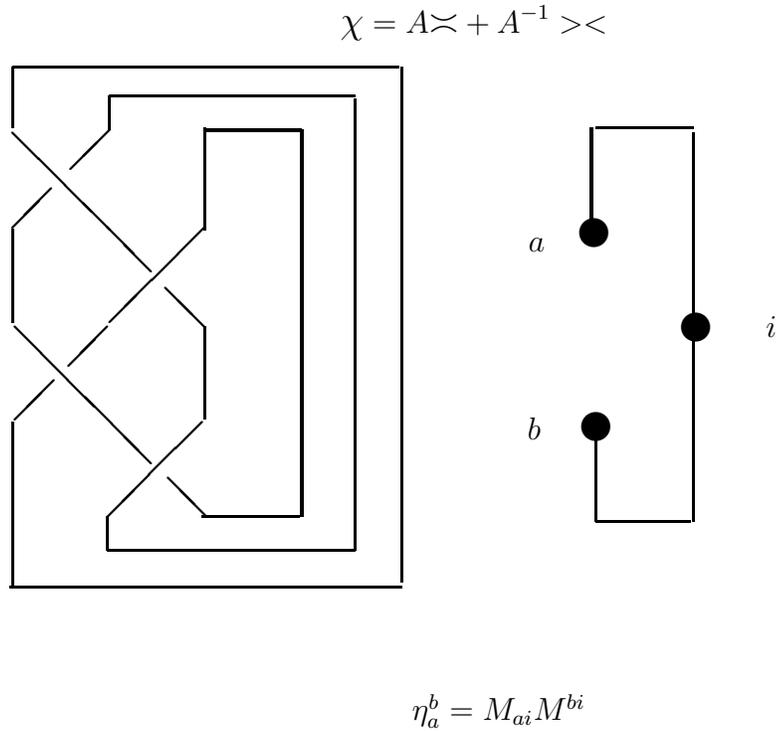
{\bf Figure 19 - Pairing Maxima and Minima for Braid Closures} 
\end{center}
\vspace{3mm}

\noindent Taken together with the loop value of  $-A^{2}  -  A^{-2}$ that is a consequence
of this matrix choice,
these equations can be regarded as a recursive algorithm for computing the
amplitude. This algorithm is the bracket state model for the (unnormalised) 
Jones polynomial  \cite{KA87}.  We have discussed this model in the previous sections.
\vspace{3mm}

The upshot of these remarks is that the bracket state summation can be reformulated as a matrix model as described in this
section.  Thus the values of the bracket polynomial can be regarded as generalized quantum amplitudes. Note also that the
model that we have described in this section can be seen as a generalization of the representation of the Temperley Lieb
algebra from Section 3 with basic matrix $M$ as above. In fact, in the case where the knot or link is a closure of a braid,
we can say even more. Suppose that $K = \overline{b}$ where $b$ is an $n$-strand braid in $B_{n}.$  Then the cups and caps can
be paired off as hown in Figure 19 so that 

$$Z(K) = \delta <K> = Trace(\eta^{\otimes n} \rho(b)).$$

\noindent Here 

$$\eta^{a}_{b} = \Sigma_{i} M_{bi}M^{ai}$$

\noindent so that $\eta = MM^{t}$ where $M^{t}$ denotes the transpose of the matrix $M$.

\noindent and $\rho:B_{n} \longrightarrow TL_{n}$ is the matrix representation of the Temperley Lieb algebra specified in
Section 3. The key to the workings of this representation of the bracket calculation for braids is the fact that for any pure
connection element $Q$ (obtained from a product of the $U_{i}'s$) in the Temperley-Lieb algebra $TL_{n}$, the evaluation

$$Trace(\eta^{\otimes n}Q)$$

\noindent is equal to $\delta^{\lambda(Q)}$ where $\lambda(Q)$ is the number of loops in the (braid) closure of the diagram
for $Q$. 
\vspace{3mm}

In general, if we have a linear function $TR:TL_{n} \longrightarrow Z[\delta]$ such that $TR(Q) = \delta^{\lambda(Q)}$ for
elements $Q$ as above, then $TR(\rho(b)) = \delta <\overline{b}>$ for any $n$-strand braid $b$.  
\vspace{3mm}

\section{Quantum Computing}

 In this paper I have concentrated on giving
a picture of the general framework of the Jones polynomial and how it is related to a
very general, in fact categorical, view of quantum mechanics. Many algorithms in
quantum topology are configured without regard to unitary evolution of the
amplitude since the constraint has been topological invariance rather than
conformation to physical reality. This gives rise to a host of problems  
of attempting to reformulate topological amplitudes
as quantum computations. A particular case in point is the bracket model for the
Jones polynomial.  It would be of great interest to see a reformulation of this
algorithm that would make it a quantum computation in the strict sense of quantum
computing. One way to think about this is to view the bracket model as a vacuum-vacuum amplitude as we have done in the last
section of this paper. Then it can be configured as a composition of operators (cups, caps and braiding).  If the braiding is
unitary. Then at least this part can be viewed as a quantum computation.
\vspace{3mm}

To see how this can be formulated consider the vacuum-vacuum computation of a link amplitude as we
have described it in section 4.  In Figure 20 we have indicated
an amplitude where the temporal decomposition consists first in a composition of cups (creations),
then braiding and then caps (annihilations). Thus we can write the amplitude in the form
$$Z_{K} = <CUP|M|CAP>$$

{\tt    \setlength{\unitlength}{0.92pt}
\begin{picture}(480,243)
\thinlines    \put(34,175){\line(1,0){202}}
              \put(33,97){\line(1,0){202}}
\thicklines   \put(71,97){\circle*{10}}
              \put(190,176){\circle*{10}}
              \put(152,174){\circle*{10}}
              \put(108,175){\circle*{10}}
              \put(150,96){\circle*{10}}
              \put(191,97){\circle*{10}}
              \put(115,98){\circle*{10}}
              \put(70,175){\circle*{10}}
              \put(149,74){\line(1,0){41}}
              \put(190,73){\line(0,1){28}}
              \put(149,98){\line(0,-1){23}}
              \put(115,77){\line(0,1){22}}
              \put(70,76){\line(1,0){45}}
              \put(70,96){\line(0,-1){20}}
              \put(152,192){\line(0,-1){19}}
              \put(110,193){\line(1,0){42}}
              \put(109,175){\line(0,1){18}}
\thinlines    \put(79,1){\framebox(283,43){$Z_{K}=<CAP|M|CUP>$}}
              \put(209,120){\framebox(39,42){M}}
              \put(277,101){\framebox(151,37){Unitary Braiding}}
              \put(269,138){\framebox(166,37){Quantum  Computation}}
              \put(254,176){\makebox(200,43){$<CAP| (Detection)$}}
              \put(255,66){\makebox(224,36){$|CUP> (Preparation)$}}
              \put(1,63){\vector(1,0){261}}
              \put(16,59){\vector(0,1){183}}
\thicklines   \put(190,218){\line(-1,0){40}}
              \put(190,178){\line(0,1){40}}
              \put(69,218){\line(1,0){41}}
              \put(69,177){\line(0,1){41}}
              \put(190,138){\line(0,-1){42}}
              \put(70,136){\line(0,-1){41}}
              \put(110,218){\line(1,0){40}}
              \put(149,98){\line(-1,1){16}}
              \put(110,135){\line(1,-1){16}}
              \put(152,138){\line(-1,-1){38}}
              \put(150,177){\line(1,-1){40}}
              \put(190,178){\line(-1,-1){16}}
              \put(151,138){\line(1,1){16}}
              \put(70,136){\line(1,1){16}}
              \put(110,176){\line(-1,-1){16}}
              \put(70,175){\line(1,-1){40}}
\end{picture}}

\begin{center}
{\bf Figure 20 - $Z_{K} = <CUP|M|CAP>$ - A Knot Quantum Computer} 
\end{center}
\vspace{3mm}

\noindent where $<CUP|$ denotes the composition of cups, $M$ is the
composition of elementary braiding matrices and $|CAP>$ is the composition of caps. We then regard
$<CUP|$ as the preparation of this state and $|CAP>$ as the detection of this state. In order to
view $Z_{K}$ as a quantum computation, we need that $M$ be a unitary operator. This will be the case
if the $R$-matrices (the solutions to the Yang-Baxter equation used in the model for this
amplitude) are unitary. In this case, each $R$-matrix can be viewed as a a quantum gate (or
possibly a composition of quantum gates) and the vacuum-vacuum diagram for the knot is interpreted
as a quantum computer. This quantum computer will probabalistically compute the values of the states
in the state sum for $Z_{K}$. In order to do so, we would need to specify those observations and preparations that corresopond to
the cups and the caps in the diagram.  A more modest proposal is to regard the braiding sector of the diagram as a quantum 
computer. That braiding sector will represent a unitary evolution, and one can ask more generally what can be computed by using such
a gate derived from a braid.
\vspace{3mm}

It should be noted that because the quantum computer gives probabilistic data, it cannot compute the knot invariants exactly.
In fact, the situation is more serious than that. If we assume that the parameters in the knot invariant are complex, then the
computer will only find (probabilistically) the absolute squares of the various complex parameters. Important phase
information will be lost, and it is not obvious that topologically invariant information about the knot can be extracted.
\vspace{3mm}

\subsection{A Unitary Representation of the Three Strand Braid Group and the Corresponding Quantum Computer}

Many questions are raised by the formulation of a quantum computer associated with a given Morse
link diagram. First of all, unitary solutions to the Yang-Baxter equation (or unitary
representations of the Artin braid group) that also give link invariants are not so easy to come by.
We gave a small example of a unitary representation of the three-strand braid group in the first section of this paper.
Thus we, are prepared to look at some aspects of the computation of a knot invariant as a quantum computation.
In fact, we can use this representation to compute the Jones polynomial for closures of 3-braids, and therefore this 
representation provides a test case for the corresponding quantum computation. We now analyse this case by first making
explicit how the bracket polynomial is computed from this representation.
\vspace{3mm}

First recall that the representation depends on two matrices $U_{1}$ and $U_{2}$ with

$$U_{1} =  \left[
\begin{array}{cc}
     \delta & 0  \\
          0 & 0
\end{array}
\right] $$

\noindent and 

$$U_{2} =  \left[
\begin{array}{cc}
     \delta^{-1} & \sqrt{1-\delta^{-2}}  \\
         \sqrt{1-\delta^{-2}}  & \delta - \delta^{-1}
\end{array}
\right]. $$

\noindent The representation is given on the two braid generators by

$$\rho(\sigma_{1})= AI + A^{-1}U_{1}$$
$$\rho(\sigma_{2})= AI + A^{-1}U_{2}$$

\noindent for any $A$ with $\delta = -A^{2}-A^{-2}$, and with 
$A=e^{i\theta}$, then $\delta = -2cos(2\theta)$. We get the specific range of
angles $|\theta| \leq \pi/6$ and $|\theta - \pi| \leq \pi/6$ that give unitary representations of
the three-strand braid group. 
\vspace{3mm}

Note that $tr(U_{1})=tr(U_{2})= \delta$ while $tr(U_{1}U_{2}) = tr(U_{2}U_{1}) =1.$
If $b$ is any braid, let $I(b)$ denote the sum of the exponents in the braid word that expresses $b$.
For $b$ a three-strand braid, it follows that 
$$\rho(b) = A^{I(b)}I + \tau(b)$$
\noindent where $I$ is the $ 2 \times 2$ identity matrix and $\tau(b)$ is a sum of products in the Temperley Lieb algebra 
involving $U_{1}$ and $U_{2}.$ Since the Temperley Lieb algebra in this dimension is generated by $I$,$U_{1}$, $U_{2}$,
$U_{1}U_{2}$ and $U_{2}U_{1}$, it follows that 
$$<\overline{b}> = A^{I(b)}\delta^{2} + tr(\tau(b))$$
\noindent where $\overline{b}$ denotes the standard braid closure of $b$, and the sharp brackets denote the bracket polynomial
as described in previous sections. From this we see at once that 
$$<\overline{b}> = tr(\rho(b)) + A^{I(b)}(\delta^{2} -2).$$
\vspace{3mm}

It follows from this calculation that the question of computing the bracket polynomial for the closure of the three-strand
braid $b$ is mathematically equivalent to the problem of computing the trace of the matrix $\rho(b).$ To what extent can our 
quantum computer determine the trace of this matrix?
\vspace{3mm}

The matrix in question is a product of unitary matrices, the quantum gates that we have associated with the braids 
$\sigma_{1}$ and $\sigma_{2}.$ The entries of the matrix $\rho(b)$ are the results of preparation and detection for the 
two dimensional basis of qubits for our machine:

$$<i|\rho(b)|j>.$$

\noindent Given that the computer is prepared in $|j>$, the probability of observing it in state $|i>$ is equal
to $|<i|\rho(b)|j>|^{2}.$ Thus we can, by running the quantum computation repeatedly, estimate the absolute squares of the
entries of the matrix $\rho(b).$  This will not yield the complex phase information that is needed for either the trace of the
matrix or the absolute value of that trace. Thus we conclude that our quantum computer can compute information relating to
the braiding process, but that it cannot approximate the full value of the bracket polynomial.
\vspace{3mm}

Note that our quantum computer does indeed have the capability to detect three strand braiding, since for a braid $b$ the
matrix $\rho(b)$ can have non-trivial off-diagaonal elements. The absolute squares of these elements are approximated by
successive runs of the quantum computer. In this quantum computer, braiding corresponds to entangled quatum states and is
dectable by that token. The bracket polynomial itself depends upon subtler phase relationships and is not detectable by this
quantum computer.
\vspace{3mm}

\subsection{Comments}

These results are less than satisfying since there does not seem to be a way to calculate the entire knot
polynomial even by probabilisitic approximations. It is not clear what the practical value of such a computation will be for
understanding a given link invariant. Nevertheless, it is to be expected that a close relationship between quantum link
invariants and quantum computing will be fruitful for both fields. \vspace{3mm} 

There are other ideas in the topology that deserve comparison with the quantum
states.  For example, topological entanglement in the sense of linking and
braiding is intuitively related to the entanglement of quantum states. In our general model using a unitary representation of the
braid group, topological entanglement entails quantum entanglement. The quantum topological
states associated with the bracket polynomial would certainly figure strongly in a quantum computing model of
this algorithm.  The specific model that we have given uses only one qubit and so does not produce entanglement. There are other 
representations of the Artin Braid group that do produce quantum entanglements corresponding to topological braiding.
These phenomena will be the subject of a subsequent paper \cite{KL}.    
\vspace{3mm}

We mention one further possibility. In the paper \cite{LB} by Lidar and Biham the authors show how to simulate special
cases of the Ising model on a quantum computer. Their method is more combinatorial and less algebraic than the approach
sketched in this section using braiding. It is possible that a generalization of their approach will work for the state sum of
the bracket polynomial. This is a topic for further research.
\vspace{3mm}

\subsection{And Quantum Field Theory}
Finally, it is important to remark that there is an interpretation of the Jones polynomial in terms of quantum field theory. 
Witten \cite{WIT89}
writes down a functional integral for link invariants in a 3-manifold  M:

$$Z(M,K) = \int dAexp[(ik/4\pi)S(M,A)] tr(Pexp(\int_{K}A)).$$

\noindent Here  $M$ denotes a 3-manifold without boundary and $A$ is a gauge field (gauge connection)  defined on $M.$
The gauge field is a one-form on $M$ with values in a representation of a  Lie algebra.    $S(M,A)$ is  the
integral over $M$ of the trace of the Chern-Simons three-form $CS = AdA + (2/3)AAA.$
(The product is the wedge product of differential forms.)
\vspace{3mm}

With the standard representation of the Lie algebra of  $SU(2)$ as $2 \times 2$
complex matrices, one can see that the formalism of
$Z(S^{},K)$  ($S^{3}$ denotes the three-dimensional sphere.) yields  the Jones
polynomial with the basic properties as we have discussed.  See Witten's
paper or \cite{WIT89} or \cite{KA91},\cite{KA95}.
\vspace{3mm}

The question is: How does the quantum field theory approach to the Jones polynomial relate to quantum computing?
One way to discuss this question is to reformulate (topological) quantum field theories as state summations, as we did for
the Jones polynomial, and then proceed in a fashion analogous to our amplitudes discussion above. It is more challenging to
try to imagine reformulating quantum computing at the level of quantum field theory. If this were accomplished, the subject of
quantum computing and the Jones polynomial might well take a new road.

\section{Summary}

In relating quantum computing with knot polynomials the key themes are unitarity and measurement.
Much is now surely unforseen. For a good survey of quantum
computing we recommend  \cite{A} and \cite{Sam} and for another view of topological issues see \cite{FR98} and
\cite{FLZ}. See \cite{Sch} for an excellent treatment of measurement theory in quantum mechanics and a useage of the Dirac formalism that 
is in resonance with the concerns of this paper.
 

\vspace{3mm}



\begin{thebibliography}{99}

\bibitem{A} D. Aharonov. Quantum computation. quant-phys/9812037 15Dec 1998.

\bibitem{AT90} M.F. Atiyah. The Geometry and Physics of Knots.  Cambridge University
Press (1990).

\bibitem{BA82}  R.J. Baxter.  Exactly Solved Models in Statistical Mechanics.  Acad.
Press (1982).

\bibitem{D58}  P.A.M. Dirac. Principles of Quantum Mechanics.  Oxford University Press
(1958).

\bibitem{FEY65}  R. Feynman and A.R. Hibbs.  Quantum Mechanics and Path Integrals. 
McGraw Hill (1965).

\bibitem{FR98} M. Freedman. Topological Views on Computational Complexity.  Documenta
Mathematica - Extra Volume ICM (1998) , pp. 453-464.

\bibitem{FLZ} M. Freedman, M. Larsen and Z. Wang. A modular functor which is universal for quantum
computation.  arXiv:quant-ph/0001108v2 1 Feb 2000.

\bibitem{JO85}  V.F.R.Jones.A polynomial invariant for links via von Neumann algebras.
Bull.Amer.Math.Soc. 129 (1985) 103-112.

\bibitem{GF91} D.S.Freed and R.E.Gompf.  Computer Calculation of WittenÕs 3-Manifold
Invariant.  Commun. Math. Phys. 141 (1991), pp. 79-117.

\bibitem{JO86}  V.F.R.Jones. A new knot polynomial and von Neumann algebras. Notices of
AMS 33 (1986) 219-225.

\bibitem{KA87}  L.H.Kauffman. State Models and the Jones Polynomial. Topology 26 (1987)
395-407.

\bibitem{KAU87} L.H. Kauffman.  On Knots.  Annals of Mathematics Studies Number 115, 
Princeton University Press (1987).

\bibitem{KA88}  L.H.Kauffman. New invariants in the theory of knots. Amer. Math. Monthly
Vol.95,No.3,March 1988. pp 195-242.

\bibitem{KAU89}  L.H.Kauffman.  Statistical mechanics and the Jones polynomial.  AMS
Contemp. Math. Series  (1989), Vol. 78. pp. 263-297.

\bibitem{KA91} L.H.  Kauffman. Knots and Physics ,  World Scientific  Pub.  (1991 and
1993).

\bibitem{KA95} L. H. Kauffman.  Functional integration and the theory of knots. J. Math.
Phys. Vol. 36, No.5 (19950, pp. 2402-2429.

\bibitem{KA98} L. H. Kauffman.  WittenÕs Integral and the Kontsevich Integral. In
Particles Fields and Gravitation. edited by Jakub Rembielinski, AIP Proceedings
No. 453, American Inst. of Physics Pub. (1998), pp. 368 - 381.

\bibitem{KAU95} L.H.Kauffman (Editor).  Knots and Applications, World Scientific Pub. Co.
(1995).

\bibitem{KA96} L.H.Kauffman. Knots and Statistical Mechanics, Proceedings of Symposia in
Applied Mathematics - The Interface of Knots and Physics (edited by L. Kauffman),
Vol.51 (1996), pp. 1-87.

\bibitem{KL} L.H.Kauffman and S. Lomonaco. Topological Entanglement and Quantum Entanglement.
(in preparation)

\bibitem{LB} D. Lidar and O. Biham. Simulating ising spin glasses on a quantum computer.
quant-ph/9611038v6 23 Sept. 1997.

\bibitem{Sam} S. Lomonaco. {\em A rosetta stone for quantum mechanics with an introduction to quantum computiation}.
quant-ph/0007045 (2000). 

\bibitem{Reid} K. Reidemeister.  {\em Knotentheorie}. Chelsea, New York (1948), 
Julius Springer (1932).

\bibitem{RT91} N.Yu.Reshetikhin and V.G. Turaev.  Invariants of 3-manifolds via link
polynomials and quantum groups.  Invent. Math. Vol. 103 (1991), pp. 547-597.

\bibitem{Sch} J. Schwinger. {\em Quantum Mechanics: Symbolism of Atomic Measurement.}
Springer-Verlag, 2001.

\bibitem{WIT89}  E. Witten. Quantum field theory and the Jones polynomial.
Commun.Math.Phys.  121 , 351-399 (1989).



\end{thebibliography}
 \end{document}